\documentclass[12pt]{article}
\usepackage{amssymb, amsmath, amscd, amsthm}
\usepackage[all]{xy}

\newtheorem{thm}{Theorem}[section]
\newtheorem{lem}[thm]{Lemma}
\newtheorem{prop}[thm]{Proposition}
\newtheorem{cor}[thm]{Corollary}

\theoremstyle{definition}
\newtheorem{defn}[thm]{Definition}
\newtheorem{exmp}[thm]{Example}

\theoremstyle{remark}

\numberwithin{equation}{section}

\def\action{\, {\scriptscriptstyle \stackrel{\circ}{{}}} \, }

  \def\C{\mathbb C}  \def\R{\mathbb R}
\def\Z{\mathbb Z}    \def\bP{\mathbb P}
\def\N{\mathbb N} 

\def\a{\alpha}  \def\b{\beta}        
\def\c{\theta}  \def\l{\lambda}   \def\p{\phi}     \def\s{\sigma}
\def\t{\tau}    \def\vp{\varphi}  \def\r{\rho}     \def\y{\eta}
\def\z{\zeta}   \def\w{\omega}    
\def\x{\xi}       

\def\G{\Gamma}        

      \def\cF{\mathcal F}
\def\cG{\mathcal G}      
\def\cL{\mathcal L}   \def\cM{\mathcal M}  
\def\cO{\mathcal O}      
\def\cT{\mathcal T}     \def\cV{\mathcal V} 
\def\cW{\mathcal W} \def\cZ{\mathcal Z}  

\DeclareMathAlphabet{\mathzap}{OT1}{pzc}{m}{it}
 \def\zD{\mathzap D}   \def\zE{\mathzap E}
 
\def\zL{\mathzap L}   
\def\zd{\mathzap d} 

   \def\Gh{\mathfrak h}   \def\Gm{\mathfrak m}

\def\Gn{\mathfrak n}      

\def\Ggl{\mathfrak{gl}}  \def\Gsl{\mathfrak{sl}}  
     \def\Go{\mathfrak{o}}

         \def\SL{\operatorname{SL}}

\def\Imag{\operatorname{Im}}    

\def\RP{\R\bP}  \def\CP{\C\bP} 

\def\del{\partial}

\def\({ \left( }     \def\){ \right) }
\def\[{ \left[ }     \def\]{ \right] }
\def\^{ \wedge }
\def\<{ \left\langle } \def\>{ \right\rangle }

\newcommand{\Frac}[2]{\frac{\textstyle #1}{\textstyle #2}}
\newcommand{\pd}[1]{\Frac{\del}{\del #1}}

\newcounter{mynum2}

\begin{document}

\title{Self-dual Zollfrei conformal structures with {\Large $\a$}-surface foliation}
\author{Fuminori Nakata
\thanks{This research is supported by
the 21st century COE program at Graduate School of Mathematical
Sciences, the University of Tokyo.}}

\maketitle

\begin{abstract}
 A global twistor correspondence is established for neutral self-dual conformal 
 structures with $\a$-surface foliation when the structure is close to the 
 standard structure on $S^2\times S^2$. 
 We need to introduce some singularity for the $\a$-surface foliation 
 such that the leaves intersect on a fixed two sphere. 
 In this correspondence, we prove that a natural double fibration is induced on 
 some quotient spaces which is equal to the standard double fibration 
 for the standard Zoll projective structure. 
 We also give a local general forms of neutral self-dual metrics with 
 $\a$-surface foliation. 
\end{abstract}


\section{Introduction}

C.~LeBrun and L.~J.~Mason investigated two kinds of  
twistor-type correspondences in \cite{bib:LM02} and \cite{bib:LM05}. 
One of them is the correspondence for 
the Zoll projective structure on two-dimensional manifolds (\cite{bib:LM02}). 
A projective structure is an equivalence class of torsion-free connections 
under the projective equivalence, where two torsion-free connections are 
called projectively equivalent if they have exactly the same unparameterized 
geodesics. A projective structure is called Zoll when all the maximal geodesics are 
closed. LeBrun and Mason proved that there is a one-to-one correspondence 
between 
\begin{itemize}
 \item equivalence classes of orientable Zoll projective structures 
    $(B,[\nabla])$, and 
 \item equivalence classes of totally real embeddings 
    $\iota : \RP^2\rightarrow\CP^2$, 
\end{itemize}
when they are close to the standard structures. 
Here $B$ is identified with the moduli space of holomorphic disks on $\CP^2$ 
whose boundaries are contained in $N=\iota(\RP^2)$. 

The second twistor correspondence constructed by 
LeBrun and Mason is the one for four-dimensional manifolds equipped with 
a neutral self-dual Zollfrei conformal structure (\cite{bib:LM05}). 
An indefinite metric on a four-dimensional manifold is called neutral 
when the signature is $(++--)$, and here we consider the indefinite 
conformal structures represented by such metrics. 
For a neutral metric on a four-manifold, we can define the self-duality condition 
in a similar way in the Riemannian case 
(cf.\cite{bib:DW},\cite{bib:DWII},\cite{bib:LM05}). 
An indefinite metric is called Zollfrei when all the maximal null geodesics are closed. 
In the neutral four-dimensional case, the Zollfrei condition and 
the self-dual condition depend only on the 
conformal class (\cite{bib:LM05}).
LeBrun and Mason introduced the notion of space-time orientation for 
a four-manifold with neutral metric, 
and proved that there is a one-to-one correspondence between 
\begin{itemize}
\item equivalence classes of space-time oriented  self-dual Zollfrei 
   conformal structures $(M,[g])$, and  
\item equivalence classes of totally real embeddings 
  $\iota : \RP^3\rightarrow\CP^3$, 
\end{itemize}
when they are close to the standard structures. 
Here $M$ is identified with the moduli space of holomorphic disks on $\CP^3$ 
whose boundaries are contained in $P=\iota(\RP^3)$. 
They also proved that only $S^2\times S^2$ admits a apace-time oriented 
neutral self-dual Zollfrei conformal structure. 

These two twistor correspondences are real, non-analytic and global versions 
of two of the three twistor correspondences
explained in \cite{bib:Hitchin82} by N.~J.~Hitchin. 
The three are twistor correspondences for 
(1) complex surfaces with projective structure, 
(2) complex four-manifolds with anti-self-dual conformal structure and 
(3) complex three-manifolds with Einstein-Weyl structure. 
The corresponding twistor space is given by 
complex manifolds $Z$ with an embedded $\CP^1$ whose normal bundle is 
$\cO(1),\cO(1)\oplus\cO(1)$ or $\cO(2)$ respectively. 
The geometric structures (1), (2) and (3) 
are given as naturl structures on the moduli spaces 
of such embedded $\CP^1$ in $Z$. 
Hitchin's argument is local, and the description is based on holomorphic category. 
The twistor space for (1) is sometimes called mini-twistor space 
(cf.\cite{bib:Calderbank}). 

The twistor correspondence for (2) is originally discovered by R.~Penrose 
(\cite{bib:Penrose76}), and the Riemannian version of this 
twistor correspondence is given by  
M.~F.~Atiyah, N.~J.~Hitchin and I.~M.~Singer (\cite{bib:AHS}). 
In the Riemannian case, anti-self-dual conformal structure is automatically analytic 
since the equation is elliptic. 
Moreover the family of rational curves on twistor space forms a 
globally defined foliation, and, by this reason, 
it is streight fowerd to translate the local description to the global case. 

On the other hand, in the cases of LeBrun and Mason, 
the equations have non-analytic solutions in general, and the 
family of $\CP^1$ in the twistor space dose not form a foliation 
different from the Riemannian case.  
LeBrun and Mason overcame these difficulties by using 
two techniques, the first one is to use the family of holomorphic disks 
instead of that of $\CP^1$, 
and the second one is the setting in terms of Zollfrei condition. 
Notice that the Zollfrei condition is an open condition 
in the space of neutral self-dual metrics (cf.\cite{bib:LM05}). 

Recently, there are some development concerning to the reduction of 
the neutral self-dual conformal structures on four-manifolds $(M,[g])$. 
M.~Dunajski and S.~West (\cite{bib:DW}) proved that, 
if there is a null conformal Killing vector field on $M$, 
then there is a natural null surface foliation containing this Killing field, 
and that a natural projective structure is induced on the leaf space. 
D.~Calderbank generalized this argument and weaken the assumption, 
the weakened assumption is given as a property for a null surface foliation on $M$. 
Both arguments are local, and formulated in smooth category. 
They also studied the analytic case, then they showed that, 
under these conditions, a twistor correspondence of case (2) 
induces a twistor correspondence of case (1) as a reduction. 

It would be natural to expect a theory of reduction 
for the two global twistor correspondences of LeBrun and Mason. 
The local theory of Dunajski,West and Calderbank would suggest that 
the natural class of such theory is neutral self-dual Zollfrei conformal structures 
with closed null surface foliation. 
Even the standard example on $S^2\times S^2$, however, does 
not contained in this class. 
The purpose of this paper is to set a nice class of 
neutral self-dual Zollfrei conformal structures equipped with an 
$\a$-surface foliation with some singularity explained later. 
Then we prove that there is a one-to-one correspondence similar to the 
twistor correspondence of LeBrun and Mason, 
and that the reduction works globally. 

In our situation, the induced projective structure on the leaf space is proved to 
be the standard Zoll projective structure. 
It would be interesting problem to find some different formulations 
so that non-standard Zoll projective structures 
are induced by the reduction.  

The organization of the paper is as follows. 
In Section \ref{Sect:proj} and Section \ref{Sect:neutral}, 
we review the definitions and properties about 
projective structures and neutral self-dual conformal structures respectively. 
In particular in Section \ref{Sect:neutral}, we prepare an explicit description 
without using spinor calculus, which enables us to 
establish the general forms of a neutral self-dual metrics with $\a$-surface foliation 
in Section \ref{Sect:foli}. 
In Section \ref{Sect:basic_foli}, 
we define a notion of basic $\a$-surface foliation 
which we need to carry out the reduction. 
Calderbank defined the notion of self-dual $\a$-surface foliation. 
In Appendix 1, we show that basic is equivalent to self-dual 
under the assumption of self-duality for the metric. 
The basic foliation, however, rather fits to our description. 
By using this notion, we give a simple proof of above results of 
Dunajski and West in Appendix 2. 

We treat the global situation in Section \ref{Sect:main&mini} 
and Section \ref{Sect:proof}. 
In Section \ref{Sect:main&mini}, 
We formulate a class of neutral self-dual conformal structures with a 
suitable $\a$-surface family, and we state the main theorem 
(Theorem \ref{thm:Main_theorem}). 
There is a low dimensional or mini-twistor version of the main theorem, and 
we prove this version in the rest of Section \ref{Sect:main&mini}. 
The proof of the main theorem is presented in Section \ref{Sect:proof}. 

In this article, we follow to the LeBrun and Mason's conventions 
of orientations and the terminology of $\a$-surface and $\b$-surface. 
We assume that all the manifolds and metrics are $C^\infty$, and 
that the topology of maps between manifolds is $C^{\infty}$-topology. 

\section{Projective structure} \label{Sect:proj}

Let $B$ be an oriented two dimensional manifold, and let 
$\cW=\bP(TB\otimes\C)$ and $\cW_\R=\bP(TB)$ be the projectivizations. 
Let $p:\cW\rightarrow B$ and $p_\R:\cW_\R\rightarrow B$ be the projections. 
Then every $w\in\cW\setminus\cW_\R$ corresponds to a complex line 
$L_w\subset T_bB\otimes\C$, where $b=p(w)$. 
Since $T_bB\otimes\C=L_w\oplus\bar{L}_w$,
$w$ defines a complex structure on $T_bB$. 
Let $\cW_+^{\,\circ}$ be one of the two connected components of 
$\cW\setminus\cW_\R$ whose element defines 
an orientation preserving complex structure, 
and we put $\cW_-^{\,\circ}$ to be the other component. 
Let $\cW_\pm$ be the closures of $\cW_\pm^{\,\circ}$, 
then we have : 
$$ \cW=\cW_+^{\,\circ} \cup \cW_-^{\,\circ} \cup \cW_\R
  = \cW_+\cup \cW_-. $$

Let $V\subset B$ be a coordinate neighborhood with an oriented coordinate 
$(y^0,y^1)$. By putting $\del_i=\frac{\del}{\del y^i}$, we can trivialize 
$\cW=\bP(T_\C B)$ on $\cV$ by 
\begin{equation} \label{eq:trivialization_of_cW}
 \CP^1\times V \overset{\sim}\longrightarrow \cW|_V \ :\ 
 \([\z_0:\z_1],b\) \longmapsto 
 [\ \z_0\del_0 + \z_1\del_1 \ ]_b. 
\end{equation}
Notice that $\cW_+|_V\simeq\{(\z,b)\in\CP^1\times V:
  \Imag \z\geq 0 \ \text{or}\ \z=\infty \}$, 
where $\z=\z_1/\z_0$ is the fiber coordinate. 

Let $\nabla$ be a torsion-free connection on $B$, 
then the connection form respecting the coordinate $(y^0,y^1)$ 
is given by $\Ggl(2,\R)$-valued 1-form $\w$ : 
$$ \w=(\w^i_j), \qquad \nabla\del_j = \w^i_j\del_i. $$ 
The horizontal lift of a tangent vector $e$ on $B$ at 
$\l=\l^i\del_i \in TB$ is 
\begin{equation}
 \tilde{e}=e-\w^i_j \l^j\pd{\l^i}. 
\end{equation}
Projecting to $\bP(TB)$, the horizontal lift of $e$ on $\cW$ at 
$\z=\z_1/\z_0$ is given by 
\begin{equation} \label{eq:lift_of_direction}
 \tilde{e}=e-\(\w^1_0+\z(\w^1_1-\w^0_0)-\z^2\w^0_1 \)(e)\pd{\z}. 
\end{equation}

Now we define a rank 1 distribution $\zL_\R$ on $\cW_\R$ as the tautological lifts,  
i.e. $\zL_{\R,(x,\z)}$ is the horizontal lift of the tangent line 
$\<\del_0+\z\del_1\>$, where $x\in B$ and 
$\z\in\RP^1\cong\cW_{\R,x}$ is the local fiber coordinate. 
From (\ref{eq:lift_of_direction}), we obtain $\zL_{\R}=\<\Gn\>$ where 
\begin{equation} \label{eq:Gn}
 \Gn=\del_0+\z\del_1-\(\w^1_0+\z(\w^1_1-\w^0_0)-\z^2\w^0_1 \)
   (\del_0+\z\del_1)\pd{\z}. 
\end{equation}
We can define a complex distribution $\cL$ on $\cW_+$ by $\cL=\<\Gn\>$, 
where $\Gn$ is extended to the vector field on $\cW_+$ by the 
analytic continuation for $\z\in\CP^1$. 
By definition, we have $\cL|_{\cW_\R}=\cL_\R\otimes\C$. 
If we put $\zd=\cL+\<\frac{\del}{\del \bar{\z}}\>$, then 
$\zd$ defines an almost complex structure on $\cW_+\setminus\cW_\R$ 
since $\zd$ satisfies $T\cW_+\otimes\C = \zd \oplus \bar{\zd}$ 
on $\cW_+\setminus\cW_\R$. 
 
Torsion-free connections $\nabla$ and $\nabla'$ on $B$ are called 
projectively equivalent, if and only if they defines exactly the same 
unparameterized geodesics. 
We call a projective structure for a projectively equivalent class $[\nabla]$. 
The following proposition is prooved in \cite{bib:LM02}
\begin{prop}
 \begin{enumerate}
  \item 
  $\zL$ and $\zd$ are defined only by $[\nabla]$, 
  \item 
  $\zL$ and $\zd$ are integrable. 
 \end{enumerate}
\end{prop}

\begin{defn}
 A projective structure $(B,[\nabla])$ is called Zoll if and only if 
 all of the maximal geodesics on $B$ are closed. 
\end{defn}
\begin{thm}[LeBrun-Mason,\cite{bib:LM02}] \label{thm:LM_Zoll}
 There is a one-to-one correspondence between 
 \begin{itemize}
  \item equivalence classes of oriented Zoll projective structures $(B,[\nabla])$, and
  \item equivalence classes of totally real embeddings 
    $\iota:\RP^2\rightarrow\CP^2$, 
 \end{itemize}
 when they are close to the standard structures. 
 The correspondence is characterized by a double fibration 
 $B\overset{p}\leftarrow (\cW_+,\cW_\R)\overset{q}
   \rightarrow(\CP^2,N),$
 where $N=\iota(\RP^2)$, $p$ is the projection, and  
 $q$ is a surjection which is holomorphic on the interior of $\cW_+$.
\end{thm}

The rough sketch of the proof is the following. 
If $(B,[\nabla])$ is given, then we can construct $(\cW_+,\cW_\R)$ 
equipped with a rank 1 foliation on $\cW_\R$. Collapsing this foliation, 
we obtain the space $(\CP^2,N)$. 
Conversely, if $\iota$ is given, then there is a family of holomorphic disks in 
$\CP^2$ such that the boundaries of disks are contained in $N$ and 
that this family defines a foliation on $\CP^2\setminus N$. 
We also remark that each holomorphic disk in this family is characterized by the 
condition : the relative homology class of the disk generates 
$H_2(\CP^2,N)\cong\Z$. 
We define $B$ to be the parameter space of this family. 
Then a Zoll projective structure $[\nabla]$ on $B$ is induced so that 
each closed geodesic is written in the form $p\action q^{-1}(\z)$ for 
some $\z\in N$. 
Notice that such family of holomorphic disks is uniquely determined 
as a deformation of the standard family if $\iota$ is close enough to the 
standard embedding.

\section{Neutral metric} \label{Sect:neutral}

Let $M$ be an oriented four-dimensional manifold, 
and let $g$ be a neutral metric on $M$
where a neutral metric is an indefinite metric of split signature. 
An oriented local frame $(e_0,e_1,e_2,e_3)$ of the tangent bundle $TM$ 
is called null tetrad if and only if its metric tensor $g_{ij}=g(e_i,e_j)$ is given by 
\begin{equation} g=(g_{ij})=
 \begin{pmatrix} &&&1 \\ &&-1& \\ &-1&& \\ 1&&& \end{pmatrix}. 
\end{equation}
Notice that, if $(e_i)$ is a null tetrad, then we obtain 
$g(\l,\l)=\det \begin{pmatrix} \l^0&\l^2 \\ \l^1&\l^3 \end{pmatrix}$ for a tangent vector $\l=\sum\l^i e_i$. 
When we make use of null tetrads, the structure group of $TM$ reduces to 
the Lie group 
\begin{equation}
 SO(2,2):= \left\{ P\in SL(4,\R) : {}^tPgP=g \right\}. 
\end{equation}
$SO(2,2)$ has two connected components and 
we denote $SO_0(2,2)$ for the identity component. 
\begin{defn}[cf.\cite{bib:LM05}]
 $M$ is called space-time orientable when the structure group of $TM$ 
 reduces to $SO_0(2,2)$. 
\end{defn}

Let $\SL(2,\R)_+$ and $\SL(2,\R)_-$ be copies of $\SL(2,\R)$. 
For each $(A,B)\in\SL(2,\R)_+\times\SL(2,\R)_-$, 
the transformation
$$ \begin{pmatrix} e_0&e_2 \\ e_1&e_3 \end{pmatrix} \longmapsto 
 A \begin{pmatrix} e_0&e_2 \\ e_1&e_3 \end{pmatrix} {}^tB$$
defines an element of $SO_0(2,2)$. 
In this way, we obtain a double covering 
$\SL(2,\R)_+\times\SL(2,\R)_-\rightarrow SO_0(2,2)$. 
The corresponding Lie algebra isomorphism 
$\Go(2,2)\simeq \Gsl(2,\R)_+\oplus \Gsl(2,\R)_-$ is given by
\begin{equation} \label{eq:Lie_alg_decomp}
  \begin{pmatrix} \,a\ &\,b\, &e&0 \\ c&d&0&e \\ 
  f&0&-d&b \\ 0&f&c&-a \end{pmatrix} \longmapsto 
  \begin{pmatrix} \frac{a-d}{2} &b \\ c&\frac{d-a}{2} \end{pmatrix} \oplus
  \begin{pmatrix} \frac{a+d}{2} &e \\ f&-\frac{a+d}{2} \end{pmatrix}. 
\end{equation}
Taking $M$ smaller, we can assume that $M$ is space-time 
oriented and the structure group of $TM$ lifts of $SL(2,\R)_+\times\SL(2,\R)_-$. 
Then we obtain a decomposition $TM=S^+\otimes S^-$, and 
the Levi-Civita connection $\nabla$ on $M$ induces the connections 
$\nabla^\pm$ on $S^\pm$. 
$S^\pm$ are called the positive and negative spin bundles, and $\nabla^\pm$ 
are called spin connections. 
If we take a local null tetrad $(e_i)$ on $TM$, then 
$\nabla$ is represented by the connection form $\w$, where 
$\w$ is a $\Go(2,2)$-valued 1-form, 
and the connection forms $\w^\pm$ of $\nabla^\pm$ are 
$\Gsl(2,\R)_\pm$-valued 1-forms, which are defined as the component 
of the decomposition of $\w$ by (\ref{eq:Lie_alg_decomp}). 

There is an eigen space decomposition 
$\wedge^2=\wedge^+ \oplus \wedge^-$ 
with respect to Hodge's $*$-operator, 
where $\wedge^2=\wedge^2 TM$ and $\wedge^\pm$ 
are the eigen spaces for the eigen values $\pm1$. 
Let $\cV$ be a null 2-plane in $T_xM$ and $v_1,v_2$ be the basis of $\cV$. 
Then the bivector $v_1\wedge v_2$ belongs to 
$\wedge^+$ or $\wedge^-$. 
\begin{defn}
 Let $\cV=\<v_1,v_2\>\subset T_xM$ be a null 2-plane.  
 If $v_1\wedge v_2\in\wedge^+$ then $\cV$ is called $\a$-plane, 
 and if $v_1\wedge v_2\in\wedge^-$ then $\cV$ is called $\b$-plane. 
 Let $S\subset M$ be a totally null embedded surface. 
 If every tangent plane on $S$ at each point is $\a$-plane, then 
 $S$ is called $\a$-surface. 
 $\b$-surface is defined in the similar way. 
\end{defn}

Let $(M,g)$ be a space-time oriented neutral manifold, and 
$(e_i)$ be a null tetrad on a open set $U\subset M$. 
From now on, we denote $e_2=\p_0, e_3=\p_1$ for a later convenience. 
The following lemma is checked by a direct calculation. 
\begin{lem}
 $ \wedge^+= \< \vp_1,\vp_2,\vp_3 \>, \quad 
    \wedge^-= \< \psi_1,\psi_2,\psi_3 \>, $ where
 \begin{equation} \label{eq:frames_of_wedge^pm}
 \begin{aligned}
  \vp_1 &=e_0\wedge e_1,\quad \vp_2=\phi_0\wedge \phi_1,\quad 
   \vp_3=\frac{1}{\sqrt{2}}(e_0\wedge \phi_1 - e_1\wedge \phi_0), \\
  \psi_1 &=e_0\wedge \phi_0,\quad \psi_2=e_1\wedge \phi_1, \quad 
   \psi_3=\frac{1}{\sqrt{2}}(e_0\wedge \phi_1 + e_1\wedge \phi_0). 
 \end{aligned}
 \end{equation}
\end{lem}

$g$ induces indefinite metrics on $\wedge^\pm$ whose metric tensors 
are both given by the following matrix with respect to the frames 
(\ref{eq:frames_of_wedge^pm}) : 
\begin{equation}
 h=\begin{pmatrix} \,0\ &1&0 \\ 1&0&0 \\ 0&0&-1 \end{pmatrix}. 
\end{equation}
Let $\Gh$ be a Lie algebra defined by
$$ \Gh = \left\{X\in\Ggl(3,\R): {}^tXh+hX=0 \right\}
 =\left\{ \begin{pmatrix} a&0&c \\ 0&-a&b \\ b&c&0
   \end{pmatrix} \right\}. $$
The Levi-Civita connection $\nabla$ induces connections on $\wedge^\pm$ 
whose connection forms are represented by $\Gh$-valued 1-form 
with respect to the frames (\ref{eq:frames_of_wedge^pm}). 

We can check that the exterior product representation associated to $\wedge^-$ is 
\begin{equation} \label{eq:omega&theta}
 \r^- : 
 \begin{pmatrix} a&b&e&0 \\ c&d&0&e \\ f&0&-d&b \\ 0&f&c&-a
  \end{pmatrix} \longmapsto 
  \begin{pmatrix} a-d&0&\sqrt{2}c \\ 0&d-a&\sqrt{2}b \\ \sqrt{2}b&\sqrt{2}c&0 
  \end{pmatrix}.  
\end{equation} 
So the connection form of the induced connection on $\wedge^-$ is given by the 
$\Gh$-valued 1-form $\c=(\c^i_j)=\r^-(\w)$, where $\w$ is the connection form of 
the Levi-Civita connection. 
This connection naturally induces the connection on $\wedge^-_\C$, 
where $\wedge^-_\C=\wedge^-\otimes\C$ is the complexification. 
The horizontal lift of a 
tangent vector $e$ on $M$ at $\l^i\psi_i\in\wedge^-_\C$ is 
\begin{equation} \label{eq:lift_to_wedge^-}
 \tilde{e}=e- \c^i_j(e) \l^j \pd{\l^i}. 
\end{equation}

Let $\cZ=\{[\psi]\in\bP(\wedge^-_\C) : g(\psi,\psi)=0 \}$ and
$\cZ_\R=\{[\psi]\in\bP(\wedge^-) : g(\psi,\psi)=0 \}$. 
Let $p:\cZ\rightarrow M$ and $p_\R:\cZ_\R\rightarrow M$ be the projections. 
Then a trivialization of $\cZ$ on the open set $U\subset M$ is given by : 
\begin{equation} \label{eq:trivialization_of_cZ}
 \iota : \CP^1\times U \overset{\sim}\longrightarrow \cZ|_U \ :\  
 \([\z_0:\z_1],x\) \longmapsto 
 [\ \z_0^2 \psi_1 + \z_1^2 \psi_2 +\sqrt{2}\z_0\z_1\psi_3\ ]_x. 
\end{equation}
This is nothing but the correspondence between 
the fiber coordinate $[\z_0:\z_1]\in\CP^1$ and the complex $\b$-plane 
$\< \z_0e_0+\z_1e_1,\ \z_0\p_0+\z_1\p_1 \>, $ 
since we have
$$ (\z_0e_0+\z_1e_1)\wedge(\z_0\p_0+\z_1\p_1)=
 \z_0^2 \psi_1 + \z_1^2 \psi_2 +\sqrt{2}\z_0\z_1\psi_3. $$
Restricting the fiber coordinate $[\z_0:\z_1]$ to $\RP^1$, 
we also obtain a trivialization of $\cZ_\R$, and each point in $\cZ_\R$  
corresponds to a real $\b$-plane in the same manner. 

Let $\b_u\subset T_xM\otimes\C$ be the complex $\b$-plane
corresponding to $u\in\cZ\setminus\cZ_\R$, where $x=p(u)$. 
Since $T_xM\otimes\C=\b_u\oplus \bar{\b}_u$, 
$z$ defines a complex structure $J$ on $T_xM$, and it is easy to check that 
$J$ preserves the metric $g$. 
Let $\cZ_+^{\,\circ}$ be one of the two connected components of 
$\cZ\setminus\cZ_\R$ whose element defines an orientation preserving 
complex structure, and we put $\cZ_-^{\,\circ}$ to be the other component. 
Let $\cZ_\pm$ be the closures of $\cZ_\pm^{\,\circ}$, 
then we have : 
$$ \cZ=\cZ_+^{\,\circ} \cup \cZ_-^{\,\circ} \cup \cZ_\R= \cZ_+\cup \cZ_-. $$

Let $\wedge^-_\C \overset{\pi}\rightarrow \bP(\wedge^-_\C)$ 
be the projectivization, 
then we obtain at $(\l^1,\l^2,\l^3)=(\z_0^2,\z_1^2,\sqrt{2}\z_0\z_1)$, 
\begin{equation}
 \pi_*\(\c^i_j\l^j\pd{\l^i}\) = \(b+\z(d-a)-\z^2 c\) \iota_*\(\pd{\z}\), 
\end{equation}
where $\z=\z_1/\z_0$ is the non-homogeneous coordinate. 
From (\ref{eq:lift_to_wedge^-}), the horizontal lift of the tangent vector $e$ on $M$ 
to $\cZ$ is 
\begin{equation} \label{eq:horizontal_lift}
  \tilde{e}=e-\(b+\z(d-a)-\z^2 c\)(e) \pd{\z}. 
\end{equation}

We can define a rank 2 distribution $\zE_\R$ on $\cZ_\R$ as the tautological lifts, 
i.e. $\zE_{\R,(x,\z)}$ is the horizontal lift of the $\b$-plane 
$\<e_0+\z e_1, \p_0+\z \p_1\>$, where $x\in M$ and 
$\z\in\RP^1\cong\cZ_{\R,x}$. $\zE_\R$ is called the twistor distribution
\cite{bib:DW} or the Lax distribution \cite{bib:Calderbank}.
From (\ref{eq:horizontal_lift}), we obtain $\zE_\R=\< \Gm_1,\Gm_2 \>$ where 
\begin{equation} \label{eq:m_1&m_2}
 \begin{aligned}
  \Gm_1 &= e_0+\z e_1 +Q_1(\z)\del_\z, \quad \, 
    Q_1(\z) =-(b+\z(d-a)-\z^2 c)(e_0+\z e_1),  \\
  \Gm_2 &= \p_0+\z \p_1 + Q_2(\z)\del_\z, \quad 
    Q_2(\z) =-(b+\z(d-a)-\z^2 c)(\p_0+\z \p_1). 
 \end{aligned}
\end{equation}
We can define a complex distribution $\zE$ on $\cZ_+$ by 
$\zE=\< \Gm_1,\Gm_2 \>$, where 
$\Gm_1$ and $\Gm_2$ are extended to the vector fields on $\cZ_+$ 
analytically in $\z\in\CP^1$. 
By definition, we have $\zE|_{\cZ_\R}=\zE_\R\otimes\C$. 
If we put $\zD=\zE+\<\frac{\del}{\del\bar{\z}}\>$, then 
$\zD$ defines an almost complex structure on $\cZ_+\setminus\cZ_\R$
since $T\cZ_+\otimes\C = \zD \oplus \bar{\zD}$ 
on $\cZ_+\setminus\cZ_\R$.
The following theorem is basic and proved in \cite{bib:LM05}, and see 
also \cite{bib:DW}. 

\begin{thm}
 \begin{enumerate}
  \item $\zE$ and $\zD$ are defined only by the conformal class $[g]$, 
  \item $\zE_\R$ is Frobenius integrable if and only if $[g]$ is self-dual. 
   Moreover, the almost complex structure on $\cZ_+\setminus\cZ_\R$
   defined from $\zD$ is integrable 
   if and only if $[g]$ is self-dual. 
 \end{enumerate}
\end{thm}

\begin{defn}
 Let $(M,[g])$ be a neutral seld-dual conformal structure, then 
 $(M,[g])$ is called Zollfrei if and only if all of the maximal null geodesics on $M$ 
 are closed.  
\end{defn}

\begin{thm}[LeBrun-Mason,\cite{bib:LM05}] \label{thm:LM_Zollfrei}
 There is a one-to-one correspondence between 
 \begin{itemize}
  \item equivalence classes of 
    space-time oriented self-dual Zollfrei conformal structures $(M,[g])$, and
  \item equivalence classes of 
    totally real embeddings $\iota:\RP^3\rightarrow\CP^3$, 
 \end{itemize}
 when they are close to the standard structures.
 The correspondence is characterized by a double fibration 
 $M\overset{p}\leftarrow (\cZ_+,\cZ_\R)\overset{q}\rightarrow (\CP^3,P), $
 where $P=\iota(\RP^3)$, 
 $p$ is the projection, and $q$ is a surjection which is holomorphic 
 on the interior of $\cZ_+$. 
\end{thm}

The proof is conceptually similar to Theorem \ref{thm:LM_Zoll}. 
$M$ is defined from $\iota$ as the parameter space of the family of 
holomorphic disks in $(\CP^3,P)$ foliating $\CP^3\setminus P$. 
Such family is uniquely determined if $\iota$ is close enough to the standard 
embedding. 

\section{$\a$-surface foliation} \label{Sect:foli}

Let $(M,g)$ be a four-manifold with a neutral metric $g$, 
and let $B$ is a two-manifold. Now we study about an $\a$-surface foliation 
$\varpi : M\rightarrow B$, i.e. 
$\varpi$ is a smooth map such that each fiber $\varpi^{-1}(b)$ on $b\in B$ 
is an $\a$-surface. 

Let $x\in M$ and $b=\varpi(x)\in B$, 
then we can take a local coordinate $(x^0,x^1,x^2,x^3)$ around $x$ 
and a coordinate $(y^0,y^1)$ around $b$ so that 
$(x^0,x^1,x^2,x^3)\overset{\varpi}\mapsto (y^0,y^1)=(x^0,x^1)$. 
Let $\cV=\< \frac{\del}{\del x^2}, \frac{\del}{\del x^3} \>$ 
be the vertical vector field, and 
we use the notation $\del_{x^i}=\frac{\del}{\del x^i}$ and so on. 
\begin{prop} \label{prop:null_tetrad}
 There is a null tetrad $(e_0,e_1,\phi_0,\phi_1)$ on $TM$ 
 which satisfies  
 \begin{enumerate}
  \item $e_0=\del_{x^0}+\a_0$ and $e_1=\del_{x^1}+\a_1$ 
    for some vertical vector fields $\a_0,\a_1\in \Gamma(\cV)$,
  \item $\phi_0$ and $\phi_1$ are vertical, i.e. $\phi_0,\phi_1 \in \Gamma(\cV)$. 
 \end{enumerate}
\end{prop}
\begin{proof}
 We take $e_0$ and $e_1$ as follows. 
 Let $\cV'$ be an $\a$-plane distribution which is transverse to $\cV$ 
 at everywhere, where $\cV'$ is not necessary integrable. 
 Since $TM=\cV\oplus \cV'$, the map 
 $\varpi_*: \cV' \overset{\sim} \rightarrow \varpi^*TB $ is an isomorphism, 
 and we can take $e_0,e_1\in \Gamma(\cV')$ so that $\varpi_*(e_i)=\del_{y^i}$. 
 If we put $\a_i=e_i-\del_{x^i}$, then $\a_i\in \Gamma(\cV)$, 
 and (1) holds. 
 
 Now $\phi_0$ and $\phi_1$ are uniquely determined so that (2) holds. 
 Actually, if we put $\phi_0=a \del_{x^2} + b \del_{x^3}$, then we have  
 $$ \begin{pmatrix} -1 \\ 0 \end{pmatrix} = 
    \begin{pmatrix} g(\del_{x^1},\del_{x^2}) & g(\del_{x^1},\del_{x^3}) \\ 
  g(\del_{x^0},\del_{x^2}) & g(\del_{x^0},\del_{x^3})  \end{pmatrix} 
 \begin{pmatrix} a \\ b \end{pmatrix}, $$
 from $g(e_1,\phi_0)=-1$ and $g(e_0,\phi_0)=0$. 
 If the $2\times 2$ matrix in the right hand side is not invertible, 
 then there is a pair of real numbers $(p,q)\neq(0,0)$ such that
 $g(\del_{x^0},p\del_{x^2}+q\del_{x^3})
   =g(\del_{x^1},p\del_{x^2}+q\del_{x^3})=0$, 
 and then $g(\del_{x^i},p\del_{x^2}+q\del_{x^3})=0$ for $i=0,1,2,3$. 
 This contracts to the non-degeneracy of $g$, so the 
 matrix is invertible, and $(a,b)$ is determined uniquely. 
 $\p_1$ is determined uniquely in the similar way. 
\end{proof}

We denote $\w$ for the connection form of the Levi-Civita connection with 
respect to the null tetrad $(e_0,e_1,\p_0,\p_1)$. Then $\w$ is a 
$\Go(2,2)$-valued 1-form, and we denote the elements in the same 
way as we used in (\ref{eq:omega&theta}). 
\begin{lem}
 The following equations hold : 
 \begin{equation} \label{eq:ideal_condition_of_V}
  \begin{aligned}
   e(\p_0)&=e(\p_1)=0, \\
   e(e_0)&=a(\p_0)=c(\p_1), \\
   e(e_1)&=b(\p_0)=d(\p_1),
  \end{aligned} \qquad
  \begin{aligned}
   a(\p_1)&=c(\p_0)=0, \\
   b(\p_1)&=d(\p_0)=0,
  \end{aligned}
 \end{equation}
 \begin{equation} \label{eq:commutators}
  \begin{array}{rcrr}
   [\p_0,\p_1]&=&(b(\p_0)+d(\p_1))\p_0 & - (a(\p_0)+c(\p_1))\p_1, \\ 
  \ [e_0,\p_0]&=&-(d(e_0)+f(\p_0))\p_0 &+ c(e_0)\p_1, \\ 
  \ [e_0,\p_1]&=&(b(e_0)-f(\p_1))\p_0 &- a(e_0)\p_1, \\
  \ [e_1,\p_0]&=&-d(e_1)\p_0 &+ (c(e_1)-f(\p_0))\p_1, \\
  \ [e_1,\p_1]&=&b(e_1)\p_0 &- (a(e_1)+f(\p_1))\p_1.
  \end{array}
 \end{equation}
\end{lem}
\begin{proof}
 Since the Levi-Civita connection $\nabla$ is torsion-free, we have
 $$ \begin{aligned}
  \ [\p_0,\p_1] &=\nabla_{\p_0}\p_1 - \nabla_{\p_1}\p_0 \\
   &= \{ e(\p_0)e_1+b(\p_0)\p_0 -a(\p_0)\p_1 \}
      - \{ e(\p_1)e_0 -d(\p_1)\p_0 +c(\p_1)\p_1 \}. \end{aligned} $$
 Since $\cV=\<\p_0,\p_1\>$ is integrable, we have $[\p_0,\p_1]\in \cV$. 
 Then we obtain 
 \begin{equation} \label{eq:integrable_condition_for_V}
  e(\p_1)=e(\p_0)=0
 \end{equation}
 and the equation for $[\p_0,\p_1]$ in (\ref{eq:commutators}). 
 By the similar calculation for $[e_i,\p_j]\in \cV$, we can check all the equations. 
\end{proof}

In the rest of this section, we assume an additional condition : 
the neutral metric $g$ is self-dual. 
Then $\Gm_1$ and $\Gm_2$ defined in 
 (\ref{eq:m_1&m_2}) satisfy the following properties. 
\begin{lem} \label{lem:Q_2=0}
 $Q_2(\z)=0$ and $(\p_0+\z \p_1)Q_1(\z)=0$. 
\end{lem}
\begin{proof}
 Since $\zE$ is integrable, 
 $[\Gm_1,\Gm_2]\subset \<\Gm_1,\Gm_2\>$. Now we have 
 $$ \begin{aligned}
  \, [\Gm_1,\Gm_2] = 
   & [e_0,\p_0]+\z^2[e_1,\p_1]+\z([e_0,\p_1]+[e_1,\p_0])
      +(e_0+\z e_1)Q_2(\z) \del_\z -Q_2(\z)e_1 \\ 
	&+Q_1(\z)\p_1 - (\p_0+\z \p_1)Q_2(\z)\del_\z 
	  +(Q_1(\z)Q'_2(\z)-Q_2(\z)Q'_1(\z))\del_\z. \end{aligned} $$
 Since $[e_0,\p_0]\in \cV$ and so on, we can write 
 $[\Gm_1,\Gm_2]=\b(\z) \Gm_2$ by using some function $\b(\z)$. 
 In the same time, we obtain the required equations. 
\end{proof}

\begin{lem} \label{lem:vanishing_of_fiber_dir}
 The following equations hold :  
 \begin{equation} \label{eq:vanishing_of_fiber_dir}
  e=0, \quad a(\p_i)=b(\p_i)=c(\p_i)=d(\p_i)=0, \quad  b=c. 
  \end{equation} 
 Especially we obtain $[\p_0,\p_1]=0$ from (\ref{eq:commutators}). 
\end{lem}
\begin{proof}
 From $Q_2(\z)=0$, we have
 $$ \begin{aligned}
   b(\p_0)&=0,\quad (d-a)(\p_0)+b(\p_1)=0, \\
   c(\p_1)&=0,\quad (a-d)(\p_1)+c(\p_0)=0. 
  \end{aligned} $$ 
 Then the first and the second equations in (\ref{eq:vanishing_of_fiber_dir}) follows 
 from these equations and (\ref{eq:ideal_condition_of_V}). 

 Now let $Q_1(\z)=q_0+q_1\z+q_2\z^2+q_3\z^3$, i.e. 
 \begin{equation} \label{eq:Q_1}
  \begin{aligned}
  q_0&=-b(e_0), \ \quad  q_1=-(d-a)(e_0)-b(e_1), \\
  q_3&=c(e_1), \qquad q_2=(a-d)(e_1)+c(e_0). 
  \end{aligned}
 \end{equation}
 We can write $[\Gm_1,\Gm_2]=\b(\z)\Gm_2$ from the proof of 
 Lemma \ref{lem:Q_2=0}, and we can put $\b(\z)=\b_0+\b_1\z+\b_2\z^2$ 
 from the relation of the degree, then from a direct calculation, we obtain 
 \begin{equation}
  \begin{array}{rcrr} 
   \ [e_0,\p_0] &=& \b_0\p_0& - q_0\ \p_1, \\[1mm]
   \ [e_0,\p_1]+[e_1,\p_0] &=& \b_1\p_0& +(\b_0-q_1) \p_1, \\[1mm]
   \ [e_1,\p_1] &=& \b_2\p_0&+(\b_1-q_2) \p_1, \\[1mm]
   \ 0 &=&& (\b_2-q_3)\p_1.
  \end{array}
 \end{equation}
 Comparing with (\ref{eq:commutators}), and using (\ref{eq:Q_1}), 
 we have $b(e_i)=c(e_i)$. 
 Since we already have $b(\p_i)=c(\p_i)$, so we obtain $b=c$. 
\end{proof}

\begin{lem}
 The following equations hold : 
  \begin{equation} \label{self-duality_for_metric}
   \begin{aligned}
    \ & \p_0b(e_0)=\p_1b(e_1)=0, \\
    \ & \p_0(a-d)(e_0)=\p_1(d-a)(e_1)=\p_0b(e_1)+\p_1b(e_0), \\
	\ & \p_0(a-d)(e_1)=\p_1(d-a)(e_0). 
   \end{aligned}
  \end{equation}
\end{lem}
\begin{proof}
 Directly deduced from (\ref{eq:Q_1}) and $(\p_0+\z \p_1)Q_1(\z)=0$. 
\end{proof}

\begin{prop} \label{prop:special_coordinate}
 Let $g$ be a neutral seld-dual metric on a four-dimensional manifold $M$ and 
 $\varpi:M\rightarrow B$ be an $\a$-surface foliation. 
 Then there is a local coordinate $(x^0,x^1,x^2,x^3)$ on $M$ so that 
 $\ker\varpi_*=\< \del_{x^2},\del_{x^3}\>$ and that 
 the metric tensor for $g$ is written in the form : 
 \begin{equation} \label{eq:metric_tensor_for_special_coordinate}
  g=(g_{ij})=\begin{pmatrix}
   p&r&0&1 \\ r&q&-1&0 \\ 0&-1&0&0 \\ 1&0&0&0 \end{pmatrix}. 
 \end{equation}
 Moreover, $p,q$ and $r$ satisfy the following equations : 
 \begin{equation} \label{eq:SD_eq_for_special_coordinate}
  \left\{ \begin{aligned}
   \del_2^2 p = \del_3^2 q &=0, \\
   \del_3^2 p + \del_2^2 q &=0, \\
   \del_2^2 r + \del_2\del_3 p 
     &= \del_3^2 r + \del_2\del_3 q=0,  
  \end{aligned} \right. 
 \end{equation}
 where $\del_i=\del_{x^i}$. 
 Conversely, for any functions $p,q$ and $r$ satisfying 
 (\ref{eq:SD_eq_for_special_coordinate}), 
 the neutral metric defined by (\ref{eq:metric_tensor_for_special_coordinate})
 is self-dual and has a natural $\a$-surface foliation.   
\end{prop}
\begin{proof}
 For given $(M,g)$ and $\varpi$, we can take a coordinate 
 $(x^0,x^1,x^2,x^3)$ and a null tetrad $(e_0,e_1,\p_0,\p_1)$ 
 on $M$ as in the Proposition \ref{prop:null_tetrad}. 
 Since we have $[\p_0,\p_1]=0$ from Lemma \ref{lem:vanishing_of_fiber_dir}, 
 we can change the coordinates $x^2,x^3$ to $w^2,w^3$ so that  
 $\p_0=\del_{w^2}, \p_1=\del_{w^3}$. So we can start from 
 $\p_0=\del_{x^2}, \p_1=\del_{x^3}$. 
 Then the metric tensor is written in the form 
 (\ref{eq:metric_tensor_for_special_coordinate}), since we have 
 \begin{equation*}
  \begin{aligned}
   g(\del_0,\del_2) &= g(e_0-\a_0, \p_0)=0, \\
   g(\del_1,\del_2) &= g(e_1-\a_1,\p_0)=-1, 
  \end{aligned}
 \end{equation*}
 and so on. 

 Now we check that the metric in the form 
 (\ref{eq:metric_tensor_for_special_coordinate}) is self-dual if and only if 
 (\ref{eq:SD_eq_for_special_coordinate}) holds. 
We take a frame on $TM$ of the form  
\begin{equation} \label{eq:null_tetrad_for_special_coordinate}
 \left\{ \begin{aligned}
  e_0 &=\del_0 +\frac{1}{2} \(r\del_2 - p\del_3 \), \\
  e_1 &=\del_1 +\frac{1}{2} \(q\del_2 - r\del_3 \),
 \end{aligned} \right. \qquad\qquad 
  \left\{ \begin{aligned}
  \p_0 &=\del_2, \\[1.5eX]
  \p_1 &=\del_3.
 \end{aligned} \right. 
\end{equation}
Then $(e_0,e_1,\p_0,\p_1)$ is a null tetrad satisfying the conditions of 
Proposition \ref{prop:null_tetrad}. 
If we calculate $[e_i,\p_j]$ and so on, and comparing with (\ref{eq:commutators}), 
then we have 
\begin{equation} \label{eq:values_for_special_coordinate}
 \begin{aligned}
  2a(e_0)&=-\del_3 p, \\ 2a(e_1)&=-\del_2 p - 2\del_3 r, 
 \end{aligned} \qquad \begin{aligned}
  2b(e_0)&=\del_2 p, \\ 2b(e_1)&=-\del_3 q, 
 \end{aligned} \qquad \begin{aligned}
  2d(e_0)&=\del_3 q + 2\del_2 r \\ 2d(e_1)&=\del_2 q.
 \end{aligned}
\end{equation}
Evaluating these equations to (\ref{self-duality_for_metric}), 
we obtain (\ref{eq:SD_eq_for_special_coordinate}). 
\end{proof}

\section{Basic foliation} \label{Sect:basic_foli}

Let $(M,g)$ be a four-manifold with a neutral metric, 
and let $\varpi: M\rightarrow B$ be an $\a$-surface foliation. 
\begin{defn}
 We define that $\varpi$ is basic if and only if the curvature $\Omega^+$ 
 of the spin connection $\nabla^+$ on $S^+$ defined by 
 (\ref{eq:Lie_alg_decomp}) is basic, i.e. 
 $i(v)\Omega^+=0$ for every  vertical vector $v\in\ker\varpi_*$. 
\end{defn}

We use the same local descriptions as Section \ref{Sect:foli}. 
Then the following lemma is proved by a direct calculation. 
\begin{lem} \label{lem:basic_condition}
 If $g$ is self-dual, then $\varpi$ is basic if and only if : 
 \begin{equation} \label{eq:basic_condition}
  \p_ib(e_k)=\p_i(a-d)(e_k)=0 \qquad \text{for} \  i,k=0,1. 
 \end{equation}
 Moreover (\ref{eq:basic_condition}) is equivalent to the equations 
 $\p_iq_j=0$ for $i=0,1$ and $j=0,1,2,3$, 
 where $Q_1(\z)=q_0+q_1\z+q_2\z^2+q_3\z^3$. 
\end{lem}

\begin{prop}[cf.\cite{bib:Calderbank}] \label{prop:induced_projective_structure}
 Let $(M,g)$ be a four-manifold with a neutral self-dual metric 
 and $\varpi:M\rightarrow B$ be an $\a$-surface foliation. 
 If $\varpi$ is basic, then there is a unique projective structure $[\nabla]$ which 
 satisfies the following condition : 
 \begin{itemize}
  \item the image of each $\b$-surface by $\varpi$ is a geodesic on $B$. 
 \end{itemize}
  Conversely, if the above condition holds for some projective structure on $B$, 
  then $\varpi$ is basic. 
\end{prop}
\begin{proof}
We take coordinate neighborhoods $U\subset M$ and $V=\varpi(U)\subset B$ 
so that the coordinates are written in the manner of Section \ref{Sect:foli}. 
Let $(e_0,e_1,\p_0,\p_1)$ be the null tetrad given by 
Proposition \ref{prop:null_tetrad}. 
Using the trivialization of $\cZ_\R$, 
$(x,\z)\in \times\RP^1\cong \cZ_\R|_U$ 
corresponds to the $\b$-surface $\b(\z)=\<e_0+\z e_1,\p_0+\z\p_1\>_x$. 
Then we have $\varpi_*(\b(\z))=\<\del_{y^0}+\z\del_{y^1}\>_{\varpi(x)}$, 
and this is a line in $T_{\varpi(x)}B$ corresponds to the point 
$(\varpi(x),\z)\in V\times\RP^1\simeq\cW_\R$. 
In this way, we obtain a map $\Pi_\R : \cZ_\R \rightarrow \cW_\R$ 
which extends holomorphically to the map 
$\Pi : \cZ_+\rightarrow\cW_+$. 

Using the above coordinates, we have 
\begin{equation} \label{eq:Pi_*(Gm_i)}
 \begin{aligned}
  \Pi_*(\Gm_1) &= \del_0+\z\del_1+Q_1(\z)\del_\z, \quad
   Q_1(\z)=q_0+q_1\z+q_2\z^2+q_3\z^3, \\
  \Pi_*(\Gm_2) &=0. 
 \end{aligned}
\end{equation}
If there is a projective structure $[\nabla]$ on $B$ satisfying 
the condition in the statement, then $\Pi_*(\zE)=\zL$, 
i.e. $\<\Gn\>=\<\Pi_*(\Gm_1)\>$. 
This equation holds only if $\p_i q_j=0$, so $\varpi$ is basic by 
Lemma \ref{lem:basic_condition}. 

Conversely, if $\varpi$ is basic, then 
$\Pi_*(\zE)$ defines a complex distribution on $\cW_+$. 
Then we can define a torsion-free connection on $B$ so that 
$\Gn=\Pi_*(\Gm_1)$. Actually, one of the examples of such connection is 
given as follows. 
From (\ref{eq:basic_condition}), $b$ and $a-d$ defines a 1-form on $V$, 
then we define the connection form $(\w^i_j)$ by 
$\w^0_1=\w^1_0=b$, and 
$$ \begin{array}{rlrl}
 \w^0_0(\del_{y^0})&=(a-d)(\del_{y^0})+b(\del_{y^1}), 
    \qquad &\w^0_0(\del_{y^1})&=b(\del_{y^0}), \\
 \w^1_1(\del_{y^0})&=b(\del_{y^1}), 
    \qquad &\w^1_1(\del_{y^1})&=(d-a)(\del_{y^1})+b(\del_{y^0}).
 \end{array} $$
Then this connection is torsion-free and the equation $\Gn=\Pi_*(\Gm_1)$ holds 
on $V=\varpi(U)$. 
This means the condition in the statement holds. 
Since the projective structure is exactly classified by the geodesics, 
such a projective structure $[\nabla]$ is uniquely defined.  
\end{proof}

\begin{cor}
 Let $(M,g)$ be a space-time oriented neutral self-dual manifold, 
 and let $\varpi:M\rightarrow B$ be an $\a$-surface foliation. 
 Then the basic condition for $\varpi$ depends only on the 
 conformal structure of $g$. 
\end{cor}
\begin{proof}
 It is obvious since the condition in the 
 Proposition \ref{prop:induced_projective_structure} 
 depends only on the conformal class $[g]$. 
\end{proof}

\begin{exmp}
Let $(x^0,x^1,x^2,x^3)$ be a coordinate on $\R^4$, and consider a metric $g$ on 
$\R^4$ whose metric tensor $g_{ij}=g(\del_{x^i},\del_{x^j})$ is 
given by : 
$$ g=(g_{ij})=\begin{pmatrix}
    p & r & 0&1 \\
	r & p & -1&0 \\
    0&-1 &0&0 \\ 1&0 &0&0 	\end{pmatrix}, 
 \quad \text{where} \quad
 \left\{ \begin{aligned}
  p&=-2x^2x^3, \\ 
  r&=(x^2)^2+(x^3)^2. 
 \end{aligned} \right.  $$
Then $g$ is neutral and self-dual, 
however the $\a$-surface foliation defined from the integrable distribution  
$\cV=\< \del_{x^2},\del_{x^3}\>$ is not $basic$. 
Actually, if we take a null tetrad in the form of 
(\ref{eq:null_tetrad_for_special_coordinate}), 
then we have
\begin{equation} \p_0b(e_1)=\p_1b(e_0)=2 \neq 0, 
\end{equation}
so (\ref{eq:basic_condition}) does not hold. 
\end{exmp}

\section{Global structure : main theorem and the mini-twistor version} 
\label{Sect:main&mini}

In this section and Section \ref{Sect:proof}, we treat the global structure. 
From now on, we call simply ``$\b$-surface" for the maximal $\b$-surface. 
Next properties are proved by LeBrun and Mason in \cite{bib:LM05}. 

\begin{prop} \label{prop:Facts_for_beta_surface}
 Let $(M,[g])$ be a space-time oriented self-dual Zollfrei conformal structure, then 
 \begin{enumerate}
   \item any two $\b$-surfaces intersect at exactly two points, 
   \item every $\b$-surface is totally geodesic embedded $S^2$, 
   \item for every $\b$-surface $\b$, the restriction of the Levi-Civita connection 
			 of $g$ to $\b$ defines a Zoll projective structure which depends only on 
			 the conformal class $[g]$, moreover this is isomorphic to the 
			 standard Zoll projective structure on $S^2$. 
 \end{enumerate}
\end{prop}

Suppose there is a closed $\a$-surface on $M$, then it satisfies 
the following lemma. 

\begin{lem} \label{lem:property_for_alpha}
 Let $(M,[g])$ be a space-time oriented self-dual Zollfrei conformal structure, 
 and let $\a$ be a closed $\a$-surface on $M$, then 
 \begin{enumerate}
   \item $\a$ is totally geodesic embedded $S^2$, 
   \item the restriction of the Levi-Civita connection of $g$ to $\a$ defines a 
             Zoll projective structure which depends only on the conformal class $[g]$,  
   \item for every $\b$-surface $\b$, the intersection of $\a$ and $\b$ 
              is either empty or $S^1$ which is the geodesic for the projective structure 
			  induced on $\a$ and $\b$. 
 \end{enumerate}
\end{lem}
\begin{proof}
 In the same way as Proposition \ref{prop:Facts_for_beta_surface}, we can check 
 that $\a$ is totally geodesic and that $[g]$ induces a projective structure on $\a$. 
 Then (1) and (2) follow from the Zollfrei condition. 
 Let $\b$ be a $\b$-surface, then $\a\cap\b$ is totally geodesic in $M$. 
 This is either empty or one-dimensional manifold, since any $\a$-plane and 
 any $\b$-plane intersect by one-dimensional subspace at a point. 
 So this is a closed geodesic on $M$. 
 Since $\a$ and $\b$ are totally geodesic, (3) holds. 
\end{proof}

We now study about self-dual Zollfrei conformal structures with 
$\a$-surface foliation. 
\begin{defn}
 We define $\cM$ to be the set of equivalence classes of the quartet
 $(M,[g],S_\infty, \cF)$ of space-time oriented self-dual Zollfrei conformal structure 
 $(M,[g])$, a $\b$-surface $S_\infty$ and a family $\cF$ of closed $\a$-surfaces 
 which satisfies the following properties : 
 (i) every $\a$-surface $\a\in\cF$ has non empty intersection with $S_\infty$, 
 (ii) $\cF$ defines a smooth foliation on $M\setminus S_\infty$. 
 An isomorphism between two quartet is a conformal isomorphism preserving 
 $S_\infty$ and $\cF$. 
\end{defn}

\begin{defn}
  We define $\bar{\cM}$ to be the set of conformal equivalence classes of 
  space-time oriented self-dual Zollfrei conformal structures $(M,[g])$. 
  Then we have a natural forgetting map $\cM\rightarrow\bar{\cM}$. 
\end{defn}

We omit to confuse a quartet $(M,[g],S_\infty,\cF)$ with its equivalence class, 
and similar to a pair $(M,[g])$. 
\begin{lem} \label{lem:antipodal_pts_of_beta}
 Let $(M,[g],S_\infty,\cF)$ be an element of $\cM$, and let $\b$ 
 be a $\b$-surface different from $S_\infty$, 
 then $\b\cap S_\infty$ is the set of antipodal points of $\b$ 
 with respect to the induced standard Zoll projective structure on $\b$. 
 Moreover, for $\a\in\cF$, $\a\cap\b$ contains $\b\cap S_\infty$ 
 if $\a\cap\b$ is not empty. 
\end{lem}
\begin{proof}
 If we take a point in $\b\setminus S_\infty$, 
 then there is a unique $\a_1\in\cF$ which contains this point. 
 If we take the other point on $\b\setminus(S_\infty\cup\a_1)$, 
 then there is a unique $\a_2\in\b$ again which contains this point. 
 Then $\a_1\cap\b$ and $\a_2\cap\b$ are different geodesics on $\b$, 
 so $\a_1\cap\a_2\cap\b$ equals to the set of antipodal points of $\b$. 
 These points belong to both $\a_1$ and $\a_2$, so it must belong to $S_\infty$. 
 On the other hand, $\b\cap S_\infty$ is just two points, so 
 $\b\cap S_\infty=\a_1\cap\a_2\cap\b$. 
 Hence $\b\cap S_\infty$ is the set of antipodal points of $\b$. 
 The latter statement is now obvious. 
\end{proof}

\begin{lem} \label{lem:cF<=>geodesic of S_infty}
 For $(M,[g],S_\infty,\cF)\in M$, 
 each $\a$-surface in $\cF$ one-to-one corresponds with a closed geodesic 
 of $S_\infty$. 
\end{lem}
\begin{proof}
 Each $\a$-surface $\a\in\cF$ determines 
 a closed geodesic $\a\cap S_\infty$ on $S_\infty$. 
 We prove that this correspondence is bijective. 
 The injectivity follows at once since $\a$-surface is totally geodesic. 
 So we check the surjectivity. 
 It is enough to show that, for each $x\in S_\infty$ and each 
 one-dimensional subspace $l\subset T_x S_\infty$, 
 there is an $\a$-surface $\a\in\cF$ such that $T_x\a\cap T_xS_\infty =l$. 
 There is a unique $\a$-plane $H\subset T_xM$ which contains $l$, 
 and we can take one dimensional subspace $l'\subset H$ different from $l$.
 Let $c$ be a closed null geodesic of $M$ which tangent to $l'$ at $x$. 
 We can take $y\in c\setminus S_\infty$ since 
 $l'$ does not tangent to $S_\infty$. 
 Then there is a unique $\a$-surface $\a\in\cF$ containing $y$, 
 and there is a unique $\b$-surface $\b$ with $c\subset \b$. 
 Since $y\in\a\cap\b$, $\a\cap\b$ is a geodesic on $\b$. 
 Since $\a\cap\b$ is a closed geodesic on $\b$ which contains $\b\cap S_\infty$ 
 and $y$ by Lemma \ref{lem:antipodal_pts_of_beta}, this is equal to $c$ . 
 Then we have $x\in \a$ and $T_x\a=H$, 
 so we have $T_x\a\cap T_xS_\infty=l$ as required. 
\end{proof}

Let $\tilde{\cG}(S_\infty)$ be the set of oriented closed geodesics on $S_\infty$. 
$\tilde{\cG}(S_\infty)$ has natural smooth structure since 
the induced projective structure on $S_\infty$ is standard. 
$\tilde{\cG}(S_\infty)$ is diffeomorphic to $S^2$, and has natural Zoll projective 
structure induced from $S_\infty$ so that a geodesic on $\tilde{\cG}(S_\infty)$
corresponds to the set of oriented geodesics on $S_\infty$ containing one 
fixed point. 

\begin{prop} \label{prop:cM=>foliation}
 Let $(M,[g],S_\infty,\cF)$ be an element of $\cM$, then there is a natural 
 identification between $\tilde{\cG}(S_\infty)$ and 
 the leaf space $B$ of the foliation on $M\setminus S_\infty$ defined by $\cF$. 
\end{prop}
\begin{proof}
 For every $\a\in\cF$, $\a$ is diffeomorphic to $S^2$, and 
 $S_\infty\cap\a=S^1$, so $(M\setminus S_\infty)\cap\a$ 
 is disjoint union of a pair of disks. 
 Hence $M\setminus S_\infty$ is foliated by such disks. 
 Each $\a\in\cF$ has natural orientation defined from the 
 space-time orientation on $M$, so each disk of the foliation is oriented. 
 Then the natural orientation is induced on the boundary of each disk. 
 In this way, the leaf space $B$ naturally corresponds to $\tilde{\cG}(S_\infty)$. 
\end{proof} 
 
\begin{prop} \label{prop:induced_Zoll_proj_str}
 For $(M,[g],S_\infty,\cF)\in\cM$, 
 the $\a$-surface foliation on $M\setminus S_\infty$ induced from $\cF$ is basic, 
 and the projective structure induced on the leaf space $B$ is isomorphic to the 
 standard Zoll projective structure. 
\end{prop}
\begin{proof}
 We already know that $B=\tilde{\cG}(S_\infty)$ has 
 the standard Zoll projective structure induced from $S_\infty$. 
 We now check that this projective structure equals to the one 
 induced from the $\a$-surface foliation. 
 Then this $\a$-surface foliation is automatically basic from 
 Proposition \ref{prop:induced_projective_structure}. 

 Let $\b$ be any $\b$-surface different from $S_\infty$. 
 It is enough to check that the set of all the leaves intersecting with $\b$ 
 corresponds to some closed geodesic on $B$ with respect to 
 the above Zoll projective structure. 
 From Lemma \ref{lem:antipodal_pts_of_beta}, 
 an $\a$-surface $\a\in\cF$ intersects with $\b$ if and only if $\a\cap S_\infty$ 
 contains the antipodal points $\b\cap S_\infty$. 
 Hence the set of $\a$-surfaces in $\cF$ intersecting with $\b$ corresponds to 
 the set of closed geodesics on $S_\infty$ containing $\b\cap S_\infty$ under the 
 correspondence of Lemma \ref{lem:cF<=>geodesic of S_infty}. 
 Such a set is a closed geodesic on $\tilde{\cG}(S_\infty)$. 
\end{proof}

Let $\RP^n\subset\CP^n$ be the standard real submanifold. 
\begin{defn}
 We define $\cT$ to be the set of equivalence classes of the pairs $(\iota,\z_0)$
 of a totally real embedding $\iota:\RP^3\rightarrow \CP^3$ 
 and a point $\z_0\in P=\iota(\RP^3)$ which satisfy : 
 \begin{itemize}
  \item $\pi(P\setminus\{\z_0\})=\RP^2$ where 
     $\pi:\CP^3\setminus\{\z_0\} \rightarrow \CP^2$ is the natural projection, 
  \item let $\CP^1_\x=\pi^{-1}(\x)\cap\{\z_0\}$ and 
     $P_\x=\CP^1_\x\cap P$ for each $\x\in\RP^2$, 
     then $(\CP^1_\x,P_\x)$ is biholomorphic to $(\CP^1,\RP^1)$, i.e. 
	 there is a biholomorphic map $\CP^1_\x\rightarrow\CP^1$ which maps  
	 $P_\x$ to $\RP^1$. 
 \end{itemize}
\end{defn}
\begin{defn}
 We define $\bar{\cT}$ to be the set of equivalence classes of totally real embeddings  
 $\iota:\RP^3\rightarrow \CP^3$. 
 Then we have a natural forgetting map $\cT\rightarrow \bar{\cT}$. 
\end{defn}
We omit to confuse a pair $(\iota,\z_0)$ or an embedding $\iota$ 
with their equivalence classes. 
Our main theorem is the following. 
We denote $f_{\cM}:\cM\rightarrow\bar{\cM}$ 
and $f_{\cT}:\cT\rightarrow\bar{\cT}$ for the forgetting maps. 
\begin{thm} \label{thm:Main_theorem}
  Let $U\subset\bar{\cM}$ and $V\subset\bar{\cT}$ be subsets containing the 
  standard elements on which the one-to-one correspondence in the sense of 
  Theorem \ref{thm:LM_Zollfrei} holds. 
  Then there is a one-to-one correspondence between 
  $f_{\cM}^{-1}(U)$ and $f_{\cT}^{-1}(V)$ 
  which satisfies the following properties : 
  if $(M,[g],S_\infty,\cF)$ corresponds to $(\iota,\z_0)$, then 
 \begin{enumerate}
  \item $(M,[g])$ corresponds to $\iota$ in the sense of 
    Theorem \ref{thm:LM_Zollfrei}, i.e. this correspondence covers the 
	correspondence between $U$ and $V$, 
  \item the standard double fibration 
    $B\leftarrow \cW_+ \rightarrow \CP^2$ is induced by using the maps
    $\varpi:M\setminus S_\infty \rightarrow B$ 
	and $\pi:\CP^3\setminus\{\z_0\}\rightarrow\CP^2$, 
	where $\varpi$ is the  $\a$-surface foliation defined from $\cF$ and 
	$\pi$ is the projection from $\z_0$. 
 \end{enumerate}
\end{thm}

Before we start to prove Theorem \ref{thm:Main_theorem}, 
we argue about a mini-twistor version in the rest of this section. 
The situation is described in the diagram (\ref{eq:mini_diagram}). 

\begin{defn}
 We define $\cM_0$ to be the set of equivalence classes 
 of triple $(S^2,[\nabla],C)$ of 
 an oriented Zoll projective structure $(S^2,[\nabla])$ and 
 a closed geodesic $C$ on $S^2$ which satisfies 
 (i) there is a smooth involution $\sigma$ on $C$, and 
 (ii) for every $x\in C$ and every closed geodesic $c$ through $x$, 
 $c$ passes through $\s(x)$. 
 We call $\sigma(x)$ the antipodal point of $x$ and 
 denote $\bar{x}$ for $\sigma(x)$. 
\end{defn}
\begin{defn}
 We define $\bar{\cM}_0$ to be the set of equivalence classes of 
 oriented Zoll projective structures $(S^2,[\nabla])$. 
 Then we have a forgetting map $\cM_0\rightarrow \bar{\cM}_0$. 
\end{defn}

\begin{defn} \label{def:cT_0}
 We define $\cT_0$ to be the set of equivalence classes of pairs 
 $(\iota,\z_0)$ of a totally real embedding $\iota:\RP^2\rightarrow \CP^2$ 
 and a point $\z_0\in N=\iota(\RP^2)$ which satisfy : 
 \begin{itemize}
  \item $\pi(N\setminus\{\z_0\})=\RP^1$ where 
    $\pi:\CP^2\setminus\{\z_0\} \rightarrow \CP^1$ is the natural projection, 
  \item let $\CP^1_\x=\pi^{-1}(\x)\cup\{\z_0\}$ and $N_\x=\CP^1_\x\cap N$, 
    then $(\CP^1_\x,P_\x)$ is biholomorphic to $(\CP^1,\RP^1)$. 
 \end{itemize}
\end{defn}
\begin{defn}
 We define $\bar{\cT}_0$ to be the set of equivalence classes of 
 the totally real embeddings $\iota:\RP^2\rightarrow \CP^2$. 
 Then we have a forgetting map $\cT_0\rightarrow \bar{\cT}_0$. 
\end{defn}

We omit to confuse $(S^2,[\nabla],C)$ with its equivalence class and so on, 
and we denote $f_{\cM_0}:\cM_0\rightarrow\bar{\cM}_0$ and 
$f_{\cT_0}:\cT_0\rightarrow\bar{\cT}_0$ for the forgetting maps. 
\begin{thm} \label{thm:mini}
 Let $U_0\subset\bar{\cM}_0$ and $V_0\subset\bar{\cT}_0$ 
 be subsets containing the standard elements on which 
 the one-to-one correspondence in the sense of Theorem \ref{thm:LM_Zoll} holds. 
 Then there is a one-to-one correspondence between 
 $f_{\cM_0}^{-1}(U_0)$ and $f_{\cT_0}^{-1}(V_0)$ 
 which covers the correspondence between 
 $U_0$ and $V_0$. 
\end{thm}
\begin{proof}
 We start from an element $(S^2,[\nabla],C)\in\cM_0$. 
 If $(S^2,[\nabla])\in U_0$, then 
 we have a double fibration $S^2 \overset{p_1}\longleftarrow (\cW_+,\cW_\R) 
 \overset{q_1}\longrightarrow (\CP^2,N)$ from Theorem \ref{thm:LM_Zoll}, 
 where $N$ is the image of the totally real embedding 
 $\iota:\RP^2\rightarrow\CP^2$. 
 We define $\z_0\in N$ to be the point corresponding to $C$, i.e. 
 $\z_0$ is the point such that $C=p_1\action q_1^{-1}(\z_0)$. 
 Let $x\in C$ be any point, and $\bar{x}$ be its antipodal point, 
 and $D_x,D_{\bar{x}}$ be the holomorphic disks on $(\CP^2,\N)$, i.e. 
 $D_x=q_1\action p_1^{-1}(x)$ and so on. 
 Notice that $\z_0\in D_x$ and $\z_0\in D_{\bar{x}}$. 
 Since each point on $\del D_x\subset N$ corresponds to some closed geodesic 
 on $S^2$ containing $x$, and since such geodesics also contain 
 $\bar{x}$ from the definition, we have $\del D_x=\del D_{\bar{x}}$. 
 Hence $l_x=D_x\cup D_{\bar{x}}$ defines a rational curve on $\CP^2$, 
 and this is proved to be a complex line. 
 Actually, let $y\in C$ be a point different from 
 $x$ and $\bar{x}$, and $l_y$ be a rational curve defined in the same way 
 as above. Then $l_x\cap l_y=\{\z_0\}$, moreover 
 $\del D_x$ and $\del D_y$ intersects transversely in $N$, 
 so $l_x$ and $l_y$ intersects only on $\z_0$ transversely. 
 Hence $l_x$ must be a complex line. 

 Let $\pi:\CP^2\setminus\{\z_0\}\rightarrow\CP^1$ be a natural projection. 
 From the above argument, we see that 
 $\pi$ maps $l_x\setminus\{\z_0\}$ to a point. 
 $N\setminus\{\z_0\}$ is foliated by lines in the form of 
 $\del D_x\setminus\{\z_0\}$, and such line one-to-one corresponds to 
 a pair of antipodal points $\{x,\bar{x}\}$ in $C$. 
 Since $\pi(N\setminus\{\z_0\})$ is the quotient space of such line foliation, 
 it is diffeomorphic to $C/\Z_2\simeq\RP^1$. 
 Since $N$ is totally real embedded $\RP^2$, it follows that 
 $\pi(N\setminus\{\z_0\})$ is also totally real submanifold in $\CP^1$. 
 Hence $(\CP^1,\pi(N\setminus\{\z_0\}))$ is biholomorphic to $(\CP^1,\RP^1)$,  
 and $(\iota,\z_0)$ defines an element of $\cT$. 

 Next we start from $(\iota,\z_0)\in\cT$. 
 If $\iota\in V_0$, then we have a double fibration
 $S^2 \overset{p_1}\longleftarrow (\cW_+,\cW_\R) 
 \overset{q_1}\longrightarrow (\CP^2,N)$. 
 We define $C=p_1\action q_1^{-1}(\z_0)$, 
 and we prove that there is a natural involution $\sigma$ on $C$. 
 
 $(\CP^1_\x,N_\x)$ consists of two holomorphic disks 
 $D_1$ and $D_2$ with $\del D_1=\del D_2=N$  
 for every $\x\in\RP^1$. Since $D_1$ and $D_2$ define 
 generators of $H_2(\CP^2,N)$, they correspond to some points 
 $x_1$ and $x_2$ in $S_\infty$ respectively. 
 Then $x_1\in C$ from $\z_0\in\del D_1$, 
 and similarly $x_2\in C$. 
 Now all the holomorphic disks containing $\z_0$ are written in the above form, 
 so $C$ equals to the union of such pairs of points. 
 We define $\sigma$ to be the involution on $C$ interchanging such two points. 
 
 It is enough to show that every closed geodesic in $S^2$ through
 $x\in C$ always pass through $\bar{x}=\sigma(x)$. 
 This is, however, obvious because each closed geodesic through 
 $x\in C$ corresponds to some point on $\del D_x=\del D_{\bar{x}}$ 
 under the double fibration, so this geodesic also pass through $\bar{x}$. 
\end{proof}

Now we explain the following diagram concerning to Theorem \ref{thm:mini}. 
\begin{equation}\label{eq:mini_diagram}
 \xymatrix{
   & \cW_+ \ar[dl]_{p_1^{\ }}  \ar[dr]^{q_1^{\ }} & \\
   S^2 & \cW_+^r \ar@{}[u]|\cup \ar[dl]_{p_1^r} \ar[d]^{\Pi} 
     \ar[dr]^{q_1^r}  & \CP^2 \\
   S^2\setminus C \ar@{}[u]|\cup \ar[d]^{\varpi} 
     & D_+ \sqcup D_-
     \ar[dl]_{p_0^{\ }} \ar[dr]^{q_0^{\ }} & \CP^2\setminus\{\z_0\}
	 \ar@{}[u]|\cup \ar[d]^{\pi} \\
   \{b_\pm\} & & \CP^1 \\ }
\end{equation}
Let $D_+$ be one of the two holomorphic disks of 
$(\CP^1,\RP^1)$ and let $D_-$ be the other one. 
Let $q_0:D_+\sqcup D_- \rightarrow \CP^1$ be the natural map. 
Let $\{b_\pm\}$ be a set consisting of two points, and 
define a map $p_0$ by $p_0(D_\pm)=b_\pm$. 
We denote $\cW_+^r=p_1^{-1}(S^2\setminus C)$, 
and let $p_1^r$ and $q_1^r$ are the restrictions of $p_1$ and $q_1$. 
Then $\Pi$ and $\varpi$ are naturally induced from $\pi$ 
so that the diagram commutes.  
Notice that $\varpi$ maps each connected component of $S^2\setminus C$ 
to one point. 
We denote $S^2_\pm=\varpi^{-1}(b_\pm)$. 

From the proof of Theorem \ref{thm:mini}, 
each point $\x\in\RP^1$ corresponds to a pair of antipodal points $\{x,\bar{x}\}$ 
of $C$ by $\CP^1_\x=D_x\cup D_{\bar{x}}$. 
Hence there is a natural isomorphism 
$i:C/\Z_2\overset{\sim}\rightarrow\RP^1$. 
On the other hand, there is a natural map $\mu:\cW_\R^r\rightarrow C/\Z_2$ 
defined as the following way. 
Each point of $\cW_\R^r$ corresponds to a pair $(x,l)$ 
of a point $x\in S^2\setminus C$ and a closed geodesic $l$ on $S^2$ 
containing $x$. Then we define $\mu(x,l)\in C/\Z_2$ to be the intersection 
$l\cap C$. 
We have $i\action\mu=q_0\action\Pi_\R$ by definition, where 
$\Pi_\R$ is the restriction of $\Pi$ to $\cW_\R^r$.  

In this way, we have checked that $\Pi:\cW_+^r\rightarrow D_+\sqcup D_-$ 
satisfies the following conditions : 
\begin{list}{($\Pi$\,\arabic{mynum2})}%
 {\usecounter{mynum2} \itemsep 0in  \leftmargin 0.4in}
 \item $\Pi$ is smooth and $\varpi\action p_1^r =p_0\action\Pi$, 
 \item there is an isomorphism $i: C/\Z_2\rightarrow \RP^1$ satisfying 
     $i\action\mu=q_0\action\Pi_\R$, 
 \item $\Pi$ is holomorphic on $\cW_+^r\setminus \cW_\R^r$.  
\end{list}

Next lemma says that such a map $\Pi$ satisfying the above conditions 
is determined uniquely up to isomorphism. 
\begin{lem} \label{lem:uniqueness_of_Pi}
 Let $\Pi$ be the map given above, and let 
 $\Pi':\cW_+^r\rightarrow D_+\sqcup D_-$ be a map which satisfies from 
 $(\Pi\,1)$ to $(\Pi\, 3)$. 
 Then there is a holomorphic automorphism $T$ on $\CP^1$ fixing $D_\pm$ 
 and satisfying $\Pi'=\tilde{T}\action\Pi$, 
 where $\tilde{T}$ is the automorphism of $D_+\sqcup D_-$ induced from $T$. 
\end{lem}
\begin{proof}
 Let $i':C/\Z_2\rightarrow\RP^2$ be the map satisfying the condition 
 $(\Pi\,2)$ for $\Pi'$, i.e. $i'\action\mu=q_0\action\Pi'_\R$. 
 Let $x\in S^2_+$ be any point, then we have 
 $$ i\action\mu_x=q_0\action\Pi_{\R,x}, \qquad 
   i'\action\mu_x=q_0\action\Pi'_{\R,x}, $$
 where $\Pi_{\R,x},\Pi'_{\R,x}$ and $\mu_x$ are restrictions of $\Pi_\R,\Pi'_\R$ 
 and $\mu$ on $\cW_{+,x}=p_1^{-1}(x)$. 
 Since $\mu_x$ is bijective, we have 
 \begin{equation} \label{eq:def_of_T}
  (q_0\action\Pi'_{\R,x})\action (q_0\action\Pi_{\R,x})^{-1}= i'\action i^{-1}. 
 \end{equation}
 The left hand side of (\ref{eq:def_of_T}) extends holomorphically to the interior 
 of $D_+$, so $i'\action i^{-1}$ extends to a holomorphic automorphism 
 on $D_+$. In the same way, if we take $x\in S_-$, we can check that 
 $i'\action i^{-1}$ 
 extends to $D_-$ holomorphically, hence there is a holomorphic automorphism $T$ 
 on $\CP^1$ fixing $D_\pm$ and satisfying $\Pi'_x=\tilde{T}\action\Pi_x$. 
 Since $T$ dose not depend on $x\in S\setminus C$, 
 this is the required automorphism. 
\end{proof}

\begin{cor} \label{cor:mini_extension_thm}
 Suppose that a given map $\Pi$ satisfies from $(\Pi\, 1)$ to $(\Pi\, 3)$, 
 then there is a unique continuous map $\pi$ which makes the diagram 
 {\rm (\ref{eq:mini_diagram})} commute. 
 Such a map $\pi$ is equivalent to the natural projection from $\z_0$. 
\end{cor}
\begin{proof}
 The map $\Pi$ satisfying the conditions from $(\Pi\, 1)$ to $(\Pi\, 3)$ is 
 essentially unique, and this is the one defined from Theorem \ref{thm:mini}. 
 So it follows that the natural projection $\pi$ is the unique map which makes the 
 diagram (\ref{eq:mini_diagram}) commute. 
\end{proof}

\section{Proof of the main theorem} \label{Sect:proof}

First we prove the following proposition. 
\begin{prop} \label{prop:cM=>cT}
 Let $(M,[g],S_\infty,\cF)\in\cM$ be an element such that 
 $(M,[g])$ is close enough to the standard Zollfrei conformal structure. 
 Then there is a unique element of $\cT$ which satisfies the conditions in 
 Theorem \ref{thm:Main_theorem}. 
\end{prop}

Let $(M,[g],S_\infty,\cF)\in\cM$ be an element as in the above statement, 
then we have a totally real embedding $\iota:\RP^3\rightarrow\CP^3$ 
corresponding to $(M,[g])$ in the sense of Theorem \ref{thm:LM_Zollfrei}. 
Let $P=\iota(\RP^3)$, then each point on $P$ corresponds to 
some $\b$-surface, so we can define $\z_0\in P$ as the point corresponding to 
$S_\infty$. 

Let $\varpi:(M\setminus S_\infty)\rightarrow B$ be the basic $\a$-surface foliation 
induced from $\cF$. 
We have the standard Zoll projective structure on $B$ by 
Proposition \ref{prop:cM=>foliation}. 
Then we have the following diagram: 
\begin{equation}\label{eq:the_diagram}
\xymatrix{
  & (\cZ_+,\cZ_\R) \ar[dl]_{p_2^{\ }}  \ar[dr]^{q_2^{\ }} & & \\
  M & (\cZ_+^r,\cZ_\R^r) \ar@{}[u]|\cup \ar[dl]_{p_2^r}
  \ar@{.>}[d]^\Pi \ar[dr]^{q_2^r}  & (\CP^3,P) & \\
  M\setminus S_\infty \ar@{}[u]|\cup \ar[d]^\varpi & (\cW_+,\cW_\R) 
    \ar[dl]_{p_1^{\ }} \ar[dr]^{q_1^{\ }} & 
	(q_2(\cZ_+^r),q_{2,\R}(\cZ_\R^r)) \ar@{}[u]|\cup 
	\ar@{.>}[d]^{\pi'} \ar@{}[r]|\subset &
    (\CP^3\setminus\{\z_0\},P\setminus\{\z_0\}) \ar@{.>}[ld]^\pi \\
  B & & (\CP^2,\RP^2) & }
\end{equation}
where $\cZ_+$ is the disk bundle over $M$ defined in the manner of 
Section \ref{Sect:neutral}, 
$\cZ_+^r=p_2^{-1}(M\setminus S_\infty)$ is its restriction and 
$\cW_+$ is the disk bundle over $B$ defined in the manner of 
Section \ref{Sect:proj}. 
Let $B\overset{p_1}\leftarrow \cW_+ \overset{q_1}\rightarrow \CP^2$ 
be the double fibration for the standard Zoll projective structure on $B$. 

Let $\zL_\R$ be the distribution on $\cW_\R$ as in Section \ref{Sect:proj}, 
and let $\zE_\R$ be the twistor distribution on $\cZ_\R$ as in 
Section \ref{Sect:neutral}. 
Then the natural map $\Pi:\cZ_+^r\rightarrow\cW_+$ is induced 
by the proof of Proposition \ref{prop:induced_projective_structure}, and 
$\Pi$ is holomorphic on $\cZ_+^r\setminus\cZ_\R^r$. 
We also have $\Pi_*(\zE_\R)=\zL_\R$ for the restriction $\Pi_\R$ of 
$\Pi$ on $\cZ_\R^r$. 
Since $q_1$ and $q_2$ are the maps which collapse the foliations defined by 
$\zL_\R$ and $\zE_\R$, 
$\Pi$ induces a continuous map $\pi':q_2(\cZ_+^r)\rightarrow\CP^2$. 
We want to prove that 
$\pi'$ smoothly extends to the natural projection 
$\pi:\CP^3\setminus\{\z_0\}\rightarrow\CP^2$, 
and that $(\CP^1_\x,P_\x)$ is biholomorphic to $(\CP^1,\RP^1)$ for each 
$\x\in\RP^2$, where $\CP^1_\x=\pi^{-1}(\x)\cup \{\z_0\}$ and 
$P_\x=\CP^1_\x\cap P$. 

We study more detail about $\pi'$. 
Let $\a$ be an $\a$-surface in $\cF$, and let $C_\a=\a\cap S_\infty$. 
If we put $\a\setminus C_\a= \a_+\sqcup\a_-$, 
then $\a_+$ and $\a_-$ are two leaves of the 
$\a$-surface foliation $\varpi: M\setminus S_\infty\rightarrow B$. 
If we put $b_\pm=\varpi(\a_\pm)$, 
then $\{b_+,b_-\}$ is the set of antipodal points on $B$ by 
Proposition \ref{prop:induced_Zoll_proj_str} and so on. 
So the corresponding holomorphic disks 
$D_{b_\pm}=q_1\action p_1^{-1}(b_\pm)$ 
have the common boundary, and 
$\CP^1_\a=D_{b_+}\cup D_{b_-}$ is complex line in $\CP^2$. 
Then we obtain the following diagram as the restriction of (\ref{eq:the_diagram}) : 
 \begin{equation}
 \xymatrix{
   & \cZ_+|_\a \ar[dl]_{p_2^\a}  \ar[dr]^{q_2^\a} \ar@{}[d]|\cup & & \\
   \a \ar@{}[d]|\cup & \cZ_+^r|_{\a\setminus C_\a} \ar[dl]
   \ar[d]^{\Pi_\a} \ar[dr]  & Q_\a \ar@{}[d]|\cup & \\
   \a\setminus C_\a \ar[d]^{\varpi_\a} & 
     D_{b_+}\sqcup D_{b_-} \ar[dl] \ar[dr] & 
	 q_2(\cZ_+^r|_{\a\setminus C_\a}) \ar[d]^{\pi'_\a} & 
	 Q_\a\setminus\{\z_0\} \ar@{}[l]|\subset \ar@{.>}[dl]^{\pi_\a} \\
   \{b_\pm\} & & \CP^1_\a & }
 \end{equation}
where
$ \cZ_+|_\a=p_2^{-1}(\a),
  \cZ_+^r|_{\a\setminus C_\a}=p_2^{-1}(\a\setminus C\a), 
  Q_\a=q_2\action p_2^{-1}(\a), $ 
and so on. 

Since $C_\a$ is a closed geodesic on $S_\infty$ with respect to 
the standard Zoll projective structure, 
$C_\a$ has a natural involution which is the restriction of the 
involution on $S_\infty$ exchanging the antipodal points. 
Hence $(\a,[\nabla],C_\a)$ defines an element of $\cM'$, 
where $[\nabla]$ is the Zoll projective structure on $\a$ defined by 
Lemma \ref{lem:property_for_alpha}. 

\begin{lem} \label{lem:finding_mini_twistor_correspondence}
 (i) $\cZ_+^\circ |_\a=(\cZ_+|_\a)\setminus (\cZ_\R|_\a)$ is a 
  complex submanifold of $\cZ_+^\circ =\cZ_+\setminus\cZ_\R$. 
 (ii) The double fibration $\a \leftarrow \cZ_+|_\a \rightarrow Q_\a$ 
  equals to the double fibration for the Zoll projective structure on $\a$ 
  given by Theorem \ref{thm:LM_Zoll}. 
  Consequently, $Q_\a$ is biholomorphic to $\CP^2$. 
\end{lem} 
\begin{proof}
 Let $\cW^\a=\bP(T\a\otimes\C)$, and we define $\cW_\pm^\a$ as in  
 the manner of Section \ref{Sect:proj}, where 
 $\cW^\a=\cW_+^\a\cup\cW_-^\a$. 
 First we construct a diffeomorphism 
 $\r:\cW_+^\a\overset\sim\rightarrow \cZ_+|_\a$. 
 Let $x\in\a$ be any point, and take a null tetrad $\{e_0,e_1,\p_0,\p_1\}$ 
 on a open neighborhood $U\subset M$ of $x$ so that $T\a=\<\p_0,\p_1\>$. 
 We can define diffeomorphisms $\cW^\a|_{U\cap\a} \overset\sim\leftarrow 
  \CP^1\times (U\cap\a) \overset\sim\rightarrow \cZ|_{U\cap\a}$ 
 by using the trivializations of $\cW^\a$ given by (\ref{eq:trivialization_of_cW}) 
 and of $\cZ$ given by (\ref{eq:trivialization_of_cZ}). 
 In other words, this map is characterized as the correspondence 
 between a complex tangent line $l$ of $\a$ and a complex $\b$-plane $\b$ 
 so that $l\subset \b$, i.e. 
 $$ \<\p_0+\z\p_1\> \longleftrightarrow \<e_0+\z e_1, \p_0+\z\p_1\>. $$ 
 This diffeomorphism does not depend on a choice of the null tetrad, 
 hence we obtain a global diffeomorphism 
 $\cW^\a\overset\sim\rightarrow \cZ|_\a$. 
 We define $\r$ to be the restriction of this diffeomorphism on $\cW_+^\a$. 
 
 We now check that $\cZ_+^\circ|_\a$ is a complex submanifold of $\cZ_+^\circ$. 
 The complex structure on $\cZ_+^\circ$  is defined so that 
 the (0,1)-vector space is $\zD=\<\Gm_1,\Gm_2,\bar{\del_\z}\>$, 
 where $\Gm_1$ and $\Gm_2$ are the horizontal lifts of 
 $e_0+\z e_1$ and $\p_0+\z\p_1$ respectively. 
 On the other hand, the complex structure on $\cW_+^{\a\,\circ}$  is defined 
 so that the (0,1)-vector space is $\zd=\<\Gn,\bar{\del_\z}\>$, 
 where $\Gn$ is the horizontal lift of $\p_0+\z\p_1$. 
 Then we obtain $\r_*(\Gn)=\Gm_2$ and $\r_*(\bar{\del_\z})=\bar{\del_\z}$, 
 hence $\r$ is holomorphic on the interior of $\cW_+^\a$. 
 
 By the similar argument, we can check $(\r_\R)_*(\zL_\R)=\zE_\R\cap T\a$
 for the restriction of $\r$ on $\cW_\R^\a$. 
 This means $q:\cZ_+|_\a\rightarrow Q_\a$ is the map which appears in the 
 double fibration for Zoll projective structure in the sense of Theorem 
 \ref{thm:LM_Zoll}. 
\end{proof}

\begin{lem} \label{lem:Q_is_complex_submanifold}
 $\CP^2\simeq Q_\a\subset\CP^3$ is a complex submanifold.  
\end{lem}
\begin{proof}
 Let $Q_{\a,\R}=q_2(\cZ_\R|_\a)$. 
 By Lemma \ref{lem:finding_mini_twistor_correspondence}, 
 $Q_\a\setminus Q_{\a,\R}\subset\CP^3\setminus P$ is a complex submanifold. 
 So it is enough to check that, for each point $\z\in Q_{\a,\R}$, there is an open 
 neighborhood in $Q_\a$ which is a complex submanifold of $\CP^3$. 
 
 Notice that $Q_{\a,\R}\subset P$ is a smooth real submanifold. 
 This follows from the facts that 
 $q_{2,\R}^\a : \cZ_\R|_\a\rightarrow Q_{\a,\R}$ is an $S^1$-bundle,  
 $q_{2,\R} : \cZ_\R\rightarrow P$ is an $S^2$-bundle, 
 and that each fiber of $q_{2,\R}^\a$ is contained some fiber of $q_{2,\R}$ as a 
 smooth real submanifold. 

 We want to show 
 \begin{equation} \label{eq:T_xiQ}
  T_\z Q_\a = T_\z Q_{\a,\R} \oplus J(T_\z Q_{\a,\R})
 \end{equation} 
 for each $\z\in Q_{\a,\R}$, where $J$ is the complex structure on $\CP^3$. 
 Originally $J$ is defined in the following manner 
 (cf.\cite{bib:LM05} proof of Theorem 7.3). 
 We can take a non-vanishing vector field $u$ on $\cZ_\R$ 
 which spans $\ker(p_{2,\R})_*$ at every point. 
 Moreover we can assume that $j(u)$ directs the interior of $\cZ_+$ where 
 $j$ is the fiberwise complex structure of $\cZ$ with respect to the 
 $\CP^1$-bundle $p_2$. Then $J$ is defined as the linear transform defined by 
 $J (q_{2*}(u))= q_{2*}(j(u))$. 
 Now equation  (\ref{eq:T_xiQ}) follows directly from this definition. 
\end{proof}

\begin{lem} \label{lem:Pi_a_is_suitable}
 $\Pi_\a$ satisfies the conditions from $(\Pi\, 1)$ to $(\Pi\, 3)$. 
\end{lem}
\begin{proof}
 It is obvious that $\Pi_\a$ satisfies $(\Pi\, 1)$ and $(\Pi\, 3)$, 
 so we check $(\Pi\, 2)$. 
 Let $S_\infty/\Z_2$ be the set of pairs of antipodal points on $S_\infty$, 
 and we define a bijection $I: S_\infty/\Z_2\rightarrow \RP^2$ by the following. 
 For each $[x]\in S_\infty/\Z_2$, the set of closed geodesics through $x$ defines 
 a geodesic on $B=\tilde{\cG}(S_\infty)$. 
 Then we define $I([x])\in\RP^2$ to be the point corresponding to this 
 geodesic in the double fibration $B\overset{p_1}\leftarrow(\cW_+,\cW_\R)
 \overset{q_1}\rightarrow(\CP^2,\RP^2)$. 

 Since $C_\a=\a\cap S_\infty$, we can define 
 $i_\a:C_\a/\Z_2 \rightarrow \RP^1_\a$ as the restriction of $I$. 
 Then we have $i_\a([x])\in\RP^1_\a$ from the definition, 
 and we have $i_\a\action\mu=q_1\action \Pi_{\a,\R}$. 
 Actually, for example on $\a_+$, 
 each point $z\in\cZ_\R|_{\a_+}$ corresponds to a pair $(x,c)$ of 
 a point $x\in\a_+$ and a closed geodesic $c$ on $\a$ through $x$. 
 Then we have $\mu(z)=[c\cap S_\infty]$ by definition. 
 Hence $i_\a\action\mu(z)=I([c\cap S_\infty])$. 
 On the other hand, let $\b_c$ be the unique $\b$-surface containing $c$, 
 then $\Pi_{\a,\R}(z)=\Pi_\R(z)\in\cW_\R$ is the point 
 defined by $(b_+,\varpi(\b_c))$, where $b_+=\varpi(x)$ and 
 $\varpi(\b_c)$ is a closed geodesic on $B$. 
 Hence we have $q_1\action\Pi_{\a,\R}(z)=I([\b_c\cap S_\infty])$ 
 from the meaning of the double fibration for $B$. 
 Since $c\cap S_\infty=\b_c\cap S_\infty$, 
 we obtain $i_\a\action\mu=q_1\action \Pi_{\a,\R}$. 
\end{proof}

\begin{cor} \label{cor:extension_of_pi_a}
 $\pi'_\a$ continuously and uniquely extends to 
 $\pi_\a:Q_\a^2\setminus \{\z_0\} \rightarrow \CP^1$, 
 and $\pi_\a$ is equivalent to the natural projection from $\z_0$. 
\end{cor}
\begin{proof}
 Directly follows from Corollary \ref{cor:mini_extension_thm}. 
\end{proof}

\begin{lem}
 There is a unique continuous extension 
 $\pi:\CP^3\setminus\{\z_0\}\rightarrow \CP^2$ of 
 $\pi':q_2(\cZ_+^r)\rightarrow \CP^2$. 
\end{lem}
\begin{proof}
 Since $\cW_+^r$ is dense in $\CP^3\setminus\{\z_0\}$, 
 the continuous extension is unique if it exist. 
 So we prove the existence. 
 For an element $\z\not\in q_2(\cW_+^r)$, 
 we define $\pi(\z)$ as follows. 
 There is a unique $x\in S_\infty$ such that $\z\in D_x$, 
 where $D_x=q_2\action p_2^{-1}(x)$. 
 Let $\a\in\cF$ be an $\a$-surface through $x$, then $\z\in Q_\a$ 
 and we put $\pi(\z)=\pi_\a(\z)$. 
 Then $\pi(\z)$ does not depend on the choice of $\a$, since 
 $\pi_\a(\z)=i_\a([x])=I(x)$ from the proof of Lemma \ref{lem:Pi_a_is_suitable}. 

 Now we check that the above $\pi$ is continuous. 
 First, notice that $\cup_\a Q_\a=\CP^3$. Actually, 
 for any $\z\in\CP^3$, we can take $x\in M$ so that 
 $x\in p_2\action q_2^{-1}(\z)$, and if we take any $\a\in\cF$ through $x$, 
 then we obtain $\z\in Q_\a$. 
 Since $\pi$ is continuous on each $Q_\a$, 
 $\pi$ is continuous on $\CP^3\setminus\{\z_0\}$. 
\end{proof}

\begin{lem} \label{lem:l_xi_is_line}
 For each $\x\in\CP^2$, $l_\xi=\pi^{-1}(\x)\cup \{\z_0\}$ is a complex line
 in $\CP^3$. In consequence, $\pi:\CP^3\setminus\{\z_0\}\rightarrow \CP^2$ is 
 the projection. 
\end{lem}
\begin{proof}
 For each $\x\in\CP^2$, there is at least one $\a\in\cF$ such that
 $\x\in\CP^1_\a$. 
 Since $\pi^{-1}(\x)=\pi_\a^{-1}(\x)$ from the definition of $\pi$, 
 $l_\xi=\pi^{-1}_\a(\xi)\cup\{\z_0\}$ is a complex line in 
 $Q_\a\simeq\CP^2$ by Corollary \ref{cor:extension_of_pi_a}. 
 Moreover $l_\xi$ is a rational curve in $\CP^3$ 
 by Lemma \ref{lem:Q_is_complex_submanifold}. 
 
 Let $\xi'\in \CP^1_\a$ be a point different from $\xi$, 
 then $l_{\xi'}$ is a rational curve in $\CP^3$. 
 $l_\xi$ and $l_{\xi'}$ are the complex lines  in 
 $Q_\a\simeq\CP^2$ which intersect only at $\z_0$. 
 Since $Q_\a\subset\CP^3$ is an embedding, 
 $l_\xi$ and $l_{\xi'}$ intersect only at $\z_0$ in 
 $\CP^3$, and the intersection is a node. 
 Hence $l_\x$ and $l_{\x'}$ are complex lines in $\CP^3$. 
\end{proof}

\begin{proof}[Proof of \ref{prop:cM=>cT}]
 For given $(M,[g],S_\infty,\cF)$, 
 we already have a totally real embedding $\iota: \RP^3\rightarrow\CP^3$ and 
 $\z_0\in P=\iota(\RP^3)$ which satisfies $\pi(P\setminus\{\z_0\})=\RP^2$ 
 for the standard projection $\pi:\CP^3\setminus\{\z_0\}\rightarrow \CP^2$. 
 For each $\x\in\RP^2$, we put $\{x,\bar{x}\}=I^{-1}(\x)$ which is the set of 
 antipodal points of $S_\infty$. 
 Let $D_x$ and $D_{\bar{x}}$ be the holomorphic disks corresponding to 
 $x$ and $\bar{x}$, then we have $l_\x=\CP^1_\x=D_x\cup D_{\bar{x}}$. 
 Since $P_\x=\del D_x=\del D_{\bar{x}}$, $(\CP^1_\x,P_\x)$ is 
 biholomorphic to $(\CP^1,\RP^1)$. Hence $(\iota,\z_0)$ is 
 an element of $\cT$ and 
 this satisfies the required conditions in Theorem \ref{thm:Main_theorem}. 
\end{proof}

Next we prove the opposite direction of the main theorem. 
\begin{prop} \label{prop:cT=>cM}
 Let $(\iota,\z_0)\in\cT$ be an element such that 
 $\iota$ is close enough to the standard embedding. 
 Then there is a unique element of $\cM$ which satisfies the conditions in 
 Theorem \ref{thm:Main_theorem}. 
\end{prop}

Let $(\iota,\z_0)\in\cT$ be an element as in the above statement. 
Let $P=\iota(\RP^3)$, and let 
$\pi : (\CP^3\setminus\{\z_0\},P\setminus\{\z_0\})
 \rightarrow (\CP^2,\RP^2)$ be the projection.  

By Theorem \ref{thm:LM_Zollfrei}, we have a space-time oriented 
self-dual Zollfrei conformal structure $(M,[g])$ and a double fibration 
$M\overset{p_2}\leftarrow\cZ_+\overset{q_2}\rightarrow\CP^3$ 
so that $q_{2,\R}(\cZ_\R)=P$. 
Each point $x\in M$ corresponds to the holomorphic disk 
$D_x=q_2\action p_2^{-1}(x)$ in $(\CP^3,P)$, 
and each point $\z\in P$ corresponds to the $\b$-surface 
$p_2\action q_2^{-1}(\z)$ on $M$. 
We define $S_\infty$ to be the $\b$-surface corresponding to 
the point $\z_0\in P$. 
Notice that, for each $x\in M\setminus S_\infty$, 
we obtain $\z_0\not\in D_x$, i.e. $D_x$ is a holomorphic disk in 
$(\CP^3\setminus\{\z_0\},P\setminus\{\z_0\})$. 

Let $B\overset{p_1}\leftarrow(\cW_+,\cW_\R)
 \overset{q_1}\rightarrow(\CP^2,\RP^2)$ 
be the double fibration given by Theorem \ref{thm:LM_Zoll}. 
Each point in $B$ corresponds to some holomorphic disk in $(\CP^2,\RP^2)$, 
and $B$ is equipped with the standard Zoll projective structure. 

Let $B/\Z_2$ be the set of pairs of antipodal points in $B$. 
Let $b\in B/\Z_2$ be a pair of antipodal points $\{b_+,b_-\}$, 
then the corresponding holomorphic disks $D_{b_\pm}$ 
have a common boundary, so 
$\CP^1_b=D_{b_+}\cup D_{b_-}$ is a complex line in $\CP^2$. 
If we put $Q_b=\pi^{-1}(\CP^1_b)\cup\{\z_0\}$, then 
$Q_b$ is a complex plane in $\CP^3$, since $\pi$ is the projection. 
We put $N_b=P\cap Q_b$. 

\begin{lem} \label{lem:decompodant_Q}
 $(Q_b,N_b)$ and $\z_0$ define an element of $\cT_0$. 
\end{lem}
\begin{proof}
 $N_b$ is the one point compactification of $\pi_\R^{-1}(\RP^1_b)$, 
 where $\RP^1_b=\RP^2\cap\CP^1_b$ and 
 $\pi_\R : P\setminus\{\z_0\}\rightarrow\RP^2$ is the restriction of $\pi$. 
 Since $\pi_\R$ is a non trivial $\R$-fibration, 
 $N_b$ is an embedded $\RP^2$ in $Q_b$. 
 Since $\z_0\in N_b$, and since the second condition in 
 Definition \ref{def:cT_0} obviously holds, 
 $(N_b,\z_0)$ defines an element of $\cT_0$. 
\end{proof}

From Lemma \ref{lem:decompodant_Q}, we obtain the similar diagram as 
(\ref{eq:mini_diagram}) which we denote 
\begin{equation} \label{eq:diagram_b}
 \xymatrix{
  & \cW_+(b) \ar[dl]_{p_1(b)}  \ar[dr]^{q_1(b)} \ar@{}[d]|\cup & \\
  \a(b) \ar@{}[d]|\cup & \cW_+^r(b) \ar[dl]
  \ar[d]^{\Pi_b} \ar[dr]  & Q_b \ar@{}[d]|\cup  \\
  \a(b)\setminus C_b \ar[d]^{\varpi_b} & 
    D_{b_+}\sqcup D_{b_-} \ar[dl] \ar[dr] & 
    Q_b\setminus\{\z_0\} \ar[d]^{\pi} \\
  \{b_\pm\} & & \CP^1_b }
\end{equation}

\begin{lem} \label{lem:embedding_of_a(b)}
 There is a natural injection $\a(b)\rightarrow M$. Moreover 
 there is a smooth map $\varpi : M\setminus S_\infty \rightarrow B$ 
 such that the restriction of $\varpi$ on $\a(b)$ is equal to $\varpi_b$.  
\end{lem}
\begin{proof}
 Let $p\in\a(b)$ be a point and $D_p=q_1(b)\action p_1(b)^{-1}(p)$ 
 be the corresponding holomorphic disk in $(Q_b,N_b)$. 
 Let $\cL_b=\{D_p\}_{p\in\a(\b)}$ be the family of such holomorphic disks in 
 $(Q_b,N_b)$, then $\cL_b$ foliates $Q_b\setminus N_b$. 
 We will soon show that $\cup_b \cL_b$ defines a family of 
 holomorphic disks in $(\CP^3,P)$ foliating $\CP^3\setminus P$, 
 then it follows that $\a(b)$ is a subset of the moduli space $M$ 
 of holomorphic disks. 
 Moreover, $\varpi$ is naturally induced as the map between 
 the set of holomorphic disks, 
 so this is smooth and $\varpi|_{\a(b)}=\varpi_b$. 

 Now we prove that $\cup_b \cL_b$ foliates $\CP^3\setminus P$. 
 For distinct points $b,b'\in B/\Z_2$, $\CP^1_b\cap\CP^1_{b'}$ consists of 
 one point $\x\in\RP^2$. 
 Then $Q_b\cap Q_b'=\pi^{-1}(\x)\cup\{\z_0\}=\CP^1_\x$. 
 If we put $P_\x=\CP^1_x\cap P$, then $(\CP^1_\x,P_\x)$ is biholomorphic to 
 $(\CP^1,\RP^1)$ by definition. So we can write $\CP^1_\x=D_1\cup D_2$, 
 where $\del D_1=\del D_2=P_\x$. 
 As in the proof of Theorem \ref{thm:mini}, $D_i$ (i=1,2) are contained in $\cL_b$ 
 and $\cL_{b'}$. Hence $\cL_b\cup\cL_{b'}$ foliates 
 $(Q_b\cup Q_{b'})\setminus(N_b\cup N_{b'})$. 
 Since $\cup_b Q_b=\CP^3$, It follows that $\cup_b \cL_b$ 
 foliates $\CP^3\setminus P$. 
\end{proof}

It follows from Lemma \ref{lem:embedding_of_a(b)} that 
the diagram (\ref{eq:diagram_b}) is the restriction of the diagram 
(\ref{eq:the_diagram}), i.e. 
$\cW_+(b)=q_2^{-1}(Q_b)=p_2^{-1}(\a(b))$, 
$p_1(b)=p_2|_{\cW_+(b)}=q_2|_{\cW_+(b)}$ and so on. 
Now we put $\cF=\{\a(b)\}_{b\in B/\Z_2}$ which is a family of 
embedded two spheres in $M$. 
Each $\a(b)$ has a Zoll projective structure defined by 
Theorem \ref{thm:mini} and Lemma \ref{lem:decompodant_Q}. 

\begin{lem} \label{lem:conditions_of_a(b)}
 $\a(b)$ is a closed $\a$-surface. 
\end{lem}
\begin{proof}
 Each closed geodesic on $\a(b)$ is written in the form 
 $$ C(\z)=p_1(b)\action q_1(b)^{-1}(\z)$$ for some $\z\in N_b$, 
 while each $\b$-surface is written in the form 
 $\b(\z)=p_2\action q_2^{-1}(\z)$ for some $\z\in P$. 
 Hence each closed geodesic on $\a(b)$ is contained in some $\b$-surface. 
 So $\a(b)$ is totally null. 
 Since totally null surface is either $\a$-surface or $\b$-surface, 
 and since $\a(b)$ is not a $\b$-surface, this is an $\a$-surface. 
\end{proof}

\begin{proof}[Proof of \ref{prop:cT=>cM}]
 For given $(\iota,\z_0)$, we take $(M,[g],S_\infty,\cF)$ as above. 
 For each $\a(b)\in\cF$, 
 $\a(b)\cap S_\infty= C_b$ is always non-empty. 
 $\cF$ defines a smooth foliation on $M\setminus S_\infty$, 
 hence $(M,[g],S_\infty,\cF)$ is an element of $\cM$. 
 This element satisfies the conditions in Theorem \ref{thm:Main_theorem}. 
\end{proof}

Theorem \ref{thm:Main_theorem} follows 
Proposition \ref{prop:cM=>cT} and Proposition \ref{prop:cT=>cM}. 

\section{Appendix 1 : Self-dual foliation} \label{Sect:SD_foli}

In Section \ref{Sect:basic_foli}, we argued  about a basic $\a$-surface foliation, 
while Calderbank observed a self-dual $\a$-surface foliation in 
\cite{bib:Calderbank}. 
Here we check that these conditions are equivalent in the assumption of 
the self-duality condition of the metric. 

Let $(M,g)$ be a four-manifold with a neutral metric $g$ 
and let $S^-$ be the negative spin bundle, then 
an $\a$-plane distribution on $M$ one-to-one corresponds with 
a subbundle $l:L\rightarrow S^-$. 

If we fix $l:L\rightarrow S^-$ and take a connection $\nabla$ on $L$, 
then we have the covariant derivative operator 
$D^\nabla:\G(S^-\otimes L^*)\rightarrow \G(TM\otimes S^- \otimes L^*)$. 
Noticing the identification $T^*M\cong S^{+*}\otimes S^{-*}\cong S^+\otimes S^-$, 
we put $T^*M\odot S^- = S^+\otimes (S^-\odot S^-) $, 
where $S^-\odot S^-$ is the symmetric tensor. 
Composing the symmetrization to the covariant derivative $D^\nabla$, 
we obtain the twistor operator 
$\cT^\nabla:\G(S^-\otimes L^*) \rightarrow \G(TM\odot S^- \otimes L^*)$.

\begin{defn}[cf.\cite{bib:Calderbank}]
 A connection $\nabla$ on $L$ is called canonical if and only if it satisfies 
 $\cT^\nabla l=0$. 
 An  $\a$-surface foliation $\varpi$ is called self-dual if and only if, 
 for the correspnding subbundle $l:L\rightarrow S^-$, 
 (i) there is a canonical connection $\nabla$ on $L$, and (ii) $\nabla$ is self-dual. 
\end{defn}

The following property is explained in \cite{bib:Calderbank}, however we 
give the proof again by using an explicit description. 
\begin{prop} \label{prop:canonical_connection}
 Let $l:L\rightarrow S_-$ be a subbundle, 
 then the $\a$-surface distribution corresponding to $l$ is integrable 
 if and only if the canonical connection on $L$ exists. 
 The canonical connection is unique if it exists. 
\end{prop}
\begin{proof}
 Since the conditions are local at all, 
 we can assume that $S^-=M\times\R^2$ and $L=M\times \R$ 
 are trivial bundles, and that $l:L\rightarrow S^-$ is a constant section 
 $l=\begin{pmatrix} 0\\1 \end{pmatrix}\in \G(S^-\otimes L^*)$. 
 Let $\begin{pmatrix} e_0 & \p_0 \\ e_1 & \p_1 \end{pmatrix}$ be 
 a null tetrad respecting the trivialization of $S^-$, 
 then the $\a$-plane distribution corresponding to $l$ is 
 given by $\<\p_0,\p_1\>$. 
 
 We denote the Levi-Civita connection of $g$ in the same way 
 as (\ref{eq:omega&theta}). 
 Let $\nabla$ be a connection on $L$ represented by a connection 1-form $\t$, 
 then the equation $\cT^\nabla l=0$ is decomposed to the following equations : 
 \begin{equation} \label{eq:canonical_connection}
 \left\{ \begin{aligned}
   \(\frac{a+d}{2}+\t\)(e_A) &=0, \\
   e(e_A)+\(\frac{a+d}{2}+\t\)(\p_A) &=0, \qquad\qquad (A=0,1). \\
   e(\p_A) &=0,
  \end{aligned} \right. 
 \end{equation}
 So the canonical connection on $L$ exists if and only if 
 $e(\p_0)=e(\p_1)=0$. 
 This is equivalent to $[\p_0,\p_1]=0$
 as in (\ref{eq:integrable_condition_for_V}), 
 hence this holds if and only if 
 $\a$-plane distribution $\<\p_0,\p_1\>$ is integrable. 
 The uniqueness of the canonical connection is obvious from 
 (\ref{eq:canonical_connection}). 
\end{proof}

\begin{lem}
 Let $(M,g)$ be a four-manifold with a neutral metric $g$ and  
 $\varpi: M\rightarrow B$ be an $\a$-surface foliation. 
 Then $\varpi$ is self-dual if and only if the following equations hold : 
 \begin{equation} \label{eq:self-dual_foli}
  \begin{aligned}
   \p_0 (a+d)(e_0) &=\p_1(a+d)(e_1)=0, \\
   \p_0 a(e_1)+\p_1a(e_0) &=-(\p_0d(e_1)+\p_1d(e_0)). 
  \end{aligned}
 \end{equation}
\end{lem}
\begin{proof} 
 Take a null tetrad as in Proposition \ref{prop:null_tetrad}. 
 Since this null tetrad fits to the proof of 
 Proposition \ref{prop:canonical_connection}, 
 the canonical connection is defined by the 1-form $\t$ satisfying
 (\ref{eq:canonical_connection}). 
 This connection is self-dual if and only if 
 $$ d\t(e_0\wedge\p_0)=d\t(e_1\wedge\p_1)
   =d\t(e_0\wedge\p_1+e_1\wedge\p_0)=0, $$
 and it is equivalent to (\ref{eq:self-dual_foli}), since $e=0$ by 
 Lemma \ref{lem:vanishing_of_fiber_dir}. 
\end{proof}

\begin{prop}
 Let $(M,g)$ be a four-manifold with a neutral self-dual metric 
 and $\varpi: M\rightarrow B$ be an $\a$-surface foliation. 
 Then $\varpi$ is self-dual if and only if it is basic. 
\end{prop}
\begin{proof}
 We take the coordinate as in Proposition \ref{prop:special_coordinate}, 
 and denote $g$ in the form (\ref{eq:metric_tensor_for_special_coordinate}). 
 If we take a null tetrad as in (\ref{eq:null_tetrad_for_special_coordinate}), 
 then each element of the connection form is given by (\ref{eq:self-dual_foli}). 
 Noticing (\ref{eq:SD_eq_for_special_coordinate}), 
 $\varpi$ is self-dual if and only if 
 \begin{equation} \label{eq:basic=self-dual}
  \del_2\del_3 p=\del_2\del_3 q=\del_3^2 p
   =\del_2^2 q=\del_2^2 r=\del_2^2 r=0. 
 \end{equation}
 This is equivalent to the basic condition (\ref{eq:basic_condition}). 
\end{proof}

\section{Appendix 2 : null conformal Killing vector field} \label{Sect:DW}

Here we treat the case of Dunajski and West, 
i.e. the case when there is a null conformal Killing vector field on $(M,g)$. 
Our method is a little far from the general treatment of twistor theory 
(cf.\cite{bib:PR}), but the calculations are easier. 

\begin{defn}
 Let $(M,g)$ be a four-manifold with a neutral metric, 
 then a vector field $K$ on $M$ 
 is called conformal Killing vector field 
 when there is a function $\y$ on $M$ satisfying $\cL_K(g)=\y g$.  
\end{defn}
 Notice that this condition depends only on the conformal structure on $M$. 
\begin{prop}[\cite{bib:DW}] \label{prop:integrability_of_distributions}
 Let $K$ be a null conformal Killing vector field on $(M,[g])$, then 
 there is a unique $\a$-plane distribution and a unique $\b$-plane distribution on $M$ 
 which contains $K$, and these distributions are both integrable. 
\end{prop}
\noindent
For the proof of this, see \cite{bib:DW} Lemma 1, and the Remark following it. 

Let $(M,g)$ be as above and $K$ be 
a null conformal Killing vector field on $M$. 
From Proposition \ref{prop:integrability_of_distributions}, $M$ has an 
$\a$-surface foliation. Taking $M$ smaller, we can assume 
the leaf space $B$ of this $\a$-surface foliation is two-dimensional manifold. 
Then we can take the coordinates and the null tetrad 
in the manner of Proposition \ref{prop:null_tetrad}. 
Now, since $K$ is a section of the $\a$-surface distribution, 
we can write 
\begin{equation}
 K=K^0\p_0+ K^1\p_1
\end{equation}
with some functions $K^0$ and $K^1$. 
We use the same description for the Levi-Civita connection as 
(\ref{eq:omega&theta}), and, 
for the simplicity, we write $a(e_i)=a_i$ and so on. 

\begin{lem}
 A null vector field $K=K^0\p_0+ K^1\p_1$ 
 is a conformal Killing vector field if and only if there is a function $\y$ on $M$ 
 and the following conditions hold: 
 \begin{equation} \label{eq:conf_Killing_eq}
 \left\{ \quad
  \begin{aligned}
   \p_0K^1=\p_1K^0 &=0 \\
   \p_0K^0=\p_1K^1 &=\y \\
   e_0 K^1 +b_0K^0-a_0K^1 &=0 \\
   e_1 K^0 -d_1K^0 + b_1K^1 &=0 \\
   e_0K^0-d_0K^0+b_0K^1 &=e_1K^1 + b_1K^0 - a_1K^1
  \end{aligned}
 \right. 
 \end{equation}
\end{lem}
\begin{proof}
 Direct calculation from $\cL_K g=\y g$. 
\end{proof}

\begin{lem}
 Let $(M,[g])$ be a neutral self-dual conformal structure, and 
 $K$ be a null conformal Killing vector field on $M$, then  
 \begin{equation} \label{eq:equations_from_conf_Killing}
  \p_0a_1+\p_1a_0=\p_0d_1+\p_1d_0=0, \quad
  \p_0a_0=\p_1d_1=0, \quad 
  \p_0b_1=\p_1b_0=0.
 \end{equation}
\end{lem}
\begin{proof}
 Differentiating (\ref{eq:conf_Killing_eq}), 
 and using (\ref{eq:commutators}), we have 
\begin{equation} \label{eq:differential_of_eta}
 \left\{ \quad
 \begin{aligned}
   e_0\y&= - (\p_1b_0) K^0 + (\p_1a_0)K^1, \\
   e_1\y&= \phantom{-} (\p_0d_1) K^0 - (\p_0b_1) K^1, \\
   e_0\y&= \phantom{-} (\p_0(d_0+b_1))K^0 - (\p_0a_1) K^1, \\
   e_1\y&= - (\p_1d_0) K^0 +(\p_1(a_1+b_0))K^1. 
  \end{aligned}
 \right. 
 \end{equation}
 Comparing these equations, and from (\ref{self-duality_for_metric}), 
 we obtain the first and the second equation of 
 (\ref{eq:equations_from_conf_Killing}). 
 Operating $\p_0$ to the first one of (\ref{eq:differential_of_eta}), 
 and $\p_1$ to the second, we have $(\p_1b_0)\y=(\p_0b_1)\y=0$. 
 We can assume $\y\neq 0$ by changing the metric in the conformal class
 $[g]$, so we have $\p_0b_1=\p_1b_0=0$. 
\end{proof}

\begin{thm} \label{thm:DW=>basic}
 The $\a$-plane distribution defined by 
 Proposition \ref{prop:integrability_of_distributions}
 is basic. 
\end{thm}
\begin{proof}
 The condition (\ref{eq:basic_condition}) is obtained 
 directly from (\ref{eq:equations_from_conf_Killing}).  
 Hence this distribution is basic. 
\end{proof}

\vspace{2ex}
\noindent
{\large \bf Acknowledgments} \\
The author most gratefully thanks his supervisor Mikio Furuta, 
for many advices and continuous help.  

\vspace{1ex}


\vspace{5mm}
\noindent \small
Graduate School of Mathematical Sciences, University of Tokyo, \\
3-8-1 Komaba Meguro, Tokyo 153-8914, Japan \\
{\it E-mail adress} : {\tt nakata@ms.u-tokyo.ac.jp}

\end{document}